\documentclass[a4paper,twoside]{article}
\usepackage{amsmath,amsfonts,amssymb,amsthm}
\usepackage{enumitem}
\usepackage[bookmarks]{hyperref}
\usepackage[center]{titlesec}
\usepackage{
mathrsfs,dsfont}
\usepackage{tocbibind}
\usepackage{geometry}
\numberwithin{equation}{section}

\newtheorem{thm}{Theorem}[section]
\newtheorem{cj}{Conjecture}[section]
\newtheorem{defi}[thm]{Definition}
\newtheorem{lem}[thm]{Lemma}
\newtheorem{cor}[thm]{Corollary}
\newtheorem{rem}{Remark}

\def\bb{\mathbb}
\def\ca{\mathcal}

\def\scr{\mathscr}

\def\bthm{\begin{thm}\def\ethm{\end{thm}}}

\def\blem{\begin{lem}\def\elem{\end{lem}}}

\def\bitm{\begin{itemize}} \def\eitm{\end{itemize}}
\def\benu{\begin{enumerate}} \def\eenu{\end{enumerate}}

\def\bpf{\begin{proof}}\def\epf{\end{proof}}
\def\beq{\begin{equation}}\def\eeq{\end{equation}}
\def\beqs{\begin{eqnarray}}\def\eeqs{\end{eqnarray}}
\def\beqsnl{\begin{eqnarray*}}\def\eeqsnl{\end{eqnarray*}}

\DeclareMathOperator{\re}{Re}
\DeclareMathOperator{\im}{Im}
\newcommand{\abs}[1]{\left\vert#1\right\vert}

\newcommand{\ip}[1]{\left<#1\right>}
\newcommand{\norm}[1]{\left\Vert#1\right\Vert}

\usepackage{fancyhdr}
\begin{document}
\title{\textbf{
Sharp 
Hardy-Littlewood-Sobolev Inequalities on Quaternionic Heisenberg Groups}}

\footnotetext{\emph{Date:} December 22, 2013, revised on July 12, 2014.}

\author{
Michael Christ\footnote{Michael Christ is supported in part by NSF grant DMS-0901569.}
,
Heping Liu \footnote{Heping Liu is supported by National Natural Science Foundation of China under Grant \#11371036 and
the Specialized Research Fund for the Doctoral Program of Higher Education of China under Grant \#2012000110059.}
,
An Zhang\footnote{Corresponding author. An Zhang is funded by China Scholarship Council under Grant No. 201306010009.}
}
\date{}
\maketitle

\begin{abstract}
In this paper, we got several sharp Hardy-Littlewood-Sobolev-type inequalities on quaternionic Heisenberg groups
(a general form due to Folland and Stein in \cite{fs}),
using the symmetrization-free method of Frank and Lieb introduced in \cite{fl}, where they considered
the analogues on classical Heisenberg group. First, we give the sharp Hardy-Littlewood-Sobolev inequalities, both on quaternionic
Heisenberg group and its equivalent on quaternionic sphere for exponent $\lambda\ge 4$.
The extremizer, as we guess, is almost uniquely constant function on sphere. Then their dual form,
sharp conformally-invariant Sobolev inequalities and the right endpoint limit case, Log-Sobolev inequality, are also obtained. For small exponent $\lambda<4$, constant function is proved to be a local extremizer.
The conformal symmetry of the inequalities and zero center-mass technique play a critical role in the argument.
\end{abstract}

\section{Introduction}
Sharp constants and extremizers for important inequalities, especially Sobolev-type inequalities, have been studied many years ago.
It has been a general hot topic in analysis, geometry, probability, PDE and quantum field theory. They play an important role
because they almost always contain or reveal profound geometric and probabilistic information.
Vast and rich beautiful literatures has been done in this subject.
Besides so many important results on Riemannian geometry, the sub-Riemannian world is more interesting,
but far away from being absolutely understood while some conclusive results in this different framework have been obtained recently.
Among so many interesting geometric inequalities, Sobolev-type inequalities exceptionally attract more attention of analysts.
In this paper, we will discuss the Hardy-Littlewood-Sobolev-type inequalities (HLS).

The classical HLS inequality on Euclidean spaces $\bb{R}^n$ states that Riesz potential, the convolution operator
with minus power of distance (or minus fractional derivative $(-\Delta)^{-\frac{\lambda}{2}}$) is a linear operator of strong-type $(p,q)$, i.e. it's bounded from
Lebesgue spaces $L^p$ to $L^q$, if exponents satisfy
admissible condition $\frac{1}{q}+\frac{\lambda}{n}=\frac{1}{p}$.
For symmetry, we can write it in a symmetric bilinear form
\[|\iint_{\bb{R}^n\times\bb{R}^n}\frac{\overline{f(x)}g(y)}{|x-y|^\lambda}dxdy| \lesssim
\|f\|_p \|g\|_q \]
with $0<\lambda<n$, and $1<p,q<\infty$, satisfying admissible condition $\frac{1}{p}+\frac{1}{q}+\frac{\lambda}{n}=2$.
It was found by Hardy, Littlewood and Sobolev almost a century ago in \cite{hl,hl-,sob}.
The existence of extremizer was first proved  by Lieb in \cite{lie}, combining the Riesz rearrangement inequality, extended Fatou
lemma and the standard compactness argument. Another proof of existence was given by Lion in \cite{lio},
using the now standard-called ``concentration-compactness" argument.
By breaking conformal symmetry and using symmetric decreasing rearrangement, sharp
constants for conjugate exponent $p=q=\frac{2n}{2n-\lambda}$ were also given in \cite{lie} and extremizer was proved to be almost uniquely of the form \[(1+|x|^2)^{-\frac{2n-\lambda}{2}},\] ignoring constant multiples, translations and dilations.
Actually, there exists for the inequality a large conformal symmetry (invariance) group, consisting of not only translations,
dilations, rotations and constant multiples,  but also an interesting inversion $f(x)\rightarrow |x|^{2n-\lambda}f(x|x|^{-2})$.
This large conformal symmetry group of the
inequality is the main difficulty, which tells
non-uniqueness of extremizer and easily vanish
the weak limit of extremizing sequence. However, it also give some convenience to add proper additional condition to find the explicit extremzier.
A unified competing symmetry method for existence and explicit fromula of extremizers
was given by Carlen and Loss in \cite{cl}. They constructed a special strong limit using alternatively
the conformal symmetry and rearrangements to any positive $L^p$ function,
which ingeniously balanced between the ``bad" and ``good" roles of the big symmetry group.
Other symmetric rearrangement-free methods of finding the sharp constant and extremizers can be found in the work of Frank and Lieb \cite{fl--,fl-}. The first reference used inversion-positivity to get results for partial exponent $\lambda\ge n-2$
while the second one demonstrate on Euclidean space the method they use later for Heisenberg group, which is an enhanced version of techniques first used by Chang and Yang in \cite{cy}.

The analogous HLS inequality on Heisenberg group $\bb{H}^n$, parameterized by $\bb{C}^n\times \bb{R}$ with group law $uu'=(z,t)(z',t')=(z+z',t+t'+2\im z\cdot_{\bb{C}} \overline{z'})$ with $u=(z,t)\in \bb{H}^n, u\in\bb{C}^n, t\in \bb{R}$ and $z\cdot_{\bb{C}}\overline{z'}=\sum_{j=1}^n z_j\overline{z'_j}$, originating from Folland and Stein's result of fractional integral operators contained in \cite{fs},
is given by, \[|\iint_{\bb{H}^n\times\bb{H}^n}\frac{\overline{f(u)}g(v)}{|u^{-1}v|^\lambda}dudv| \lesssim \|f\|_p \|g\|_q\]
with $u^{-1}v$ the left translation of group elements $v$ by $u, ~du, dv$ the Haar measure, $|\cdot|$ the homogeneous norm, $0<\lambda<Q=2n+2$, and $1<p,q<\infty$, satisfying admissible condition $\frac{1}{p}+\frac{1}{q}+\frac{\lambda}{Q}=2$.
Here, $Q$ is the homogeneous dimension of Heisenberg group. For the sharp problem, the first important progress is made by Jerison
and Lee in the celebrated paper \cite{jl} for the special case $p=q, \lambda=Q-2$, whose extremizer is proved to be almost uniquely
$\left((1+|z|^2)^2+|t|^2\right)^{-\frac{Q+2}{4}}$ with $u=(z,t)\in \bb{H}^n, z\in\bb{C}^n, t\in\bb{R}$. Here, by ``almost", we mean modulo constant multiples, left translations and dilations. They used the extremizers to study Cauchy-Riemann
Yamabe problem, which is about solutions to the Euler-Lagrange equation of certain equivalent Sobolev inequalities functional.
Actually in \cite{jl}, they gave, in the dual form,
sharp Sobolev-embedding inequality of exponent 2,
which in this special case,
coincides with one kind of conformally-invariant Sobolev-type inequality involving intertwining operator of order 2, in other words, involving sublaplacian operator $\ca{L}$.
General Sobolev-embedding (Folland-Stein) inequality \cite{fs} states
\[\|f\|_{p*}\lesssim \|Xf\|_p,\]
with $1<p<Q, p*=\frac{pQ}{Q-p}$ and $X$ is the horizontal distributional vector fields,
while the conformally-invariant Sobolev inequality \cite{bfm} states
\[\int_{\bb{H}^n} \bar{f}\ca{L}_df \gtrsim \|f\|_p^2,\] with $0<d<Q, p=\frac{2Q}{Q-d}$ and $\ca{L}_d$ is the general intertwining operator of complementary series representation of semisimple Lie group (here for $SU(n+1,1)$).
Then still for conjugate exponent, Branson, Fontana and Morpurgo
gave in \cite{bfm} the endpoint limit case $\lambda\rightarrow 0$, a sharp Log-HLS inequality, which states that: for any normalized nonnegtive $f$ satisfying
$\int_{\bb{H}^n} f=|\bb{S}^{2n+1}|, \int_{\bb{H}^n} f(u)\log (1+|u|^2)<\infty,$  we have
\[\iint_{\bb{H}^n\times\bb{H}^n} f(u)\log\frac{2}{|u^{-1}v|^2}g(v)dudv\le \frac{|\bb{S}^{2n+1}|}{Q}\left(\int_{\bb{H}^n} f\log f+\int_{\bb{H}^n} g\log g+2\log 2|\bb{S}^{2n+1}|\right) \]
with almost unique extremizer
$\left((1+|z|^2)^2+|t|^2\right)^{-\frac{Q}{2}}$. They also conjectured that \[\left((1+|z|^2)^2+|t|^2\right)^{-\frac{2Q-\lambda}{4}}\]
is the almost unique extremizer for general $\lambda$ and ``almost" is the same meaning as before. Actually, they stated in a dual form ``conformally-invariant Sobolev inequalities" involving intertwining operators. Inspired by the former two special cases,
Frank and Lieb \cite{fl} recently killed the conjecture of \cite{bfm} (still for conjugate exponent), by proving that
a cleverly chosen extremizer is just function above. This extremizer is similar to, but also very different from that of
Euclidean HLS inequality. The level set of Euclidean extremizer is just Euclidean sphere, however,
because of the different degrees of coordinates in Heisenberg case, the level set of extremizer is neither the level set
of homogeneous norm, nor isoperimetric surface (w.r.t the perimetric measure \cite{capogna2007introduction}). So the symmetric rearrangement method for Euclidean case won't
work again for noncommutative Heisenberg case.
The existence was proved in a dual form as the distance power kernel is positive definite, combining ``refined HLS" inequality
introduced in \cite{bgx,bg}, Fatou lemma and the compactness argument. As before in Euclidean case, there is also a large
conformal symmetry group of HLS inequality on Heisenberg group, which makes it very easy for the weak limit of an extremizing sequence to vanish. However, Frank and Lieb proved that left translations and dilations are the
only ways of losing compactness, which means, we can recover a strong limit by repulling the extremizing sequence back through them. Refined HLS is a beautiful tool to realize this idea, which was first used by G\'{e}rard \cite{g} to prove the existence of extremizer of Sobolev inequality on $\bb{R}^n$.
Another proof of existence for general exponent $(p,q)$, using concentration-compactness argument from \cite{lio}, was
given by Han in \cite{h}. In finding the explicit extremizers, boundary extension of Cayley transform, which is an isomorphism
between two models, Siegel upper domain $D$ and unit ball $B$ models, of related rank one symmetric space (solvable extension of group $\bb{H}^n$),
plays a basic role in converting the inequality from group onto sphere. First, it's easier to deal on sphere. Second, it gives
a more clear view about how to choose the special extremizer that we want.
The critical step to filter the extremizer is utilizing the method from Herch, Chang and Yang \cite{her,cy} to prove a purported inverse second-variation
inequality reaches equality only by constant function.
Actually, it's a natural idea to break the huge conformal symmetry group by restricting
extremizers to functions satisfying certain ``zero center-mass" condition as in above references.
Fortunately, by checking the eigenvalues of the quadratic forms in the inverse second-variation inequality, this smartly chosen extremizers can only be constant function
and finally, by turning back, we figure out all extremziers.

Now, as conjectured in \cite{fl}, we want to extend the results to sharp HLS inequalities for general 2-step nilponent Lie groups, especially of Heisenberg type (H-type) introduced by Kaplan in \cite{k}. We may first consider several special examples of Iwasawa type (I-type \cite{cdkr}), a subclass of H-type, which is the nilpotent part of the Iwasawa decomposition of a semisimple Lie group of rank one. It's proved that there are only four cases for groups of I-type, including commutative Euclidean space, classical Heisenberg group, quaternionic Heisenberg group and octonionic Heisenberg group.
This paper is focused on quaternionic Heisenberg group $G$, another group of I-type with 3 dimensional center, parameterized by its Lie algebra $\bb{H}^n\times \im\bb{H}$ ($\bb{H}$ is the quaternions) with corresponding induced group law
$uu'=(q,w)(q',w')=(q+q',w+w'+2\im q\cdot\overline{q'})$, where $u=(q,w)\in G, u\in\bb{H}^n, w\in \im \bb{H}$ and $q\cdot\overline{q'}=\sum_{j=1}^n q\overline{q'}$ is the quaternionic product.
Accurately, we proved  that, for $\lambda\ge 4$ and conjugate exponent $p=q=\frac{2Q}{2Q-\lambda}$,
extremizer exist and is ``almost" uniquely,
constant function for sharp HLS inequality on quaternionic sphere (Theorem \ref{t-s})
and
\[\big((1+|q|^2)^2+|w|^2\big)^{-\frac{2Q-\lambda}{4}}\] for sharp HLS inequality on quaternionic Heisenberg group with group element $u=(q,w)\in G, q\in \bb{H}^n, w\in\im\bb{H}$ and $Q$ the homogeneous dimension (Theorem \ref{t-q}).
Aready known results by us are about the dual form for $\lambda=Q-2$ in a series of papers \cite{gv,imv,ivanovoptimal2012}.
They considered such sharp embedding problem when studying the sub-Riemanian Yamabe problem, which states
\[\|f\|_{\frac{2Q}{Q-2}} \lesssim \|Xf\|_2,\]
where $X$ is the horizontal distributional gradient.
\cite{gv} proved in the unifying I-type framework that restricted extremizer
satisfying additional ``partial symmetry" for above inequality is almost uniquely $\big((1+|q|^2)^2+|w|^2\big)^{-\frac{Q-2}{4}}$, which is also used there to give all the solutions, still with the addtional symmetry, of the Yamabe equation. \cite{imv} gave all extremziers
on the seven dimensional standard quaternionic sphere, while \cite{ivanovoptimal2012} obtained results for all dimensions.
Both of them got rid of the partial symmetry conditon.
So, we also extend the result of \cite{ivanovoptimal2012}. The conformal symmetry group is again the core role in the analysis.
However, a bit of difference on quaternionic case from Heisenberg case is that the special extremizer we choose can be proved to be
constant function only for partial exponent $\lambda\ge 4$. Note the special case $\lambda=Q-2=4n+4>4$ is included. For $\lambda< 4$, it seems that the chosen second-variation
and zero-center mass condition is not enough to filter the extremizers to be constant function. So, maybe more information of
the conformal symmetry group should be excavated to strengthen the condition, while Euler-Lagrange equation can also do some favor. However, instead of global consideration, we prove that constant function is still the local extremizer. The global problem will be discussed in following series papers, where octonionic and more general H-type group will also be studied.
The dual form, conformally-invariant Sobolev inequality and endpoint limit Log-Sobolev inequality are also obtained in this script.

\emph{Arrgangement of this paper:} In section \ref{p}, we will give the preliminary introducing the quaternionic Heisenberg group,
Cayley transform, which moves the problem onto sphere, and intertwining operators both on group and sphere.
Then, in section \ref{mr}, we give our main theorems, sharp HLS inequalities. Two related Sobolev-type inequalities are also obtained. We only give their proof outline: basic dual argument and sharp HLS inequality imply the sharp conformally-invariant Sobolev inequalities (actually they are equivalent),
while the endpoint limit sharp Log-Sobolev inequality is derived from the functional limit argument. In section \ref{p-hls}, we give the proof of the two main theorems about sharp HLS
inequalities. First, we use a quaternionic edition ``refined HLS" to prove the existence. Then we want to give the explicit formula
of extremizers. We characterize the symmetry group by zero center-mass condition
in subsection \ref{zcm}, reduce the problem to a second-variation inequality in subsection \ref{ss-isv},
compute the sharp constant and explicit formula of all extremizers in subsection \ref{sc}, assuming the \textbf{claim} that constant function is the almost unique extremizer.
The last but critical step of the whole proof is to prove a bilinear inequality and the condition for equality.
We give Funk-Hecke formula and compute eigenvalues of integral operators with kernel of power of distance function in subsection \ref{fh},
and at last in subsection \ref{p-isv}, use the results to give the proof of the bilinear inequality, which is just
the inverse of the second-variation inequality obtained in subsection \ref{ss-isv}.

\section{Preliminary}\label{p}
\emph{Quaternionic Heisenberg group:}
noncommutative quaternionic ring is given by \[\bb{H}=\{q=a+bi+cj+dk~|~a,b,c,d \in \bb{R}; i^2=j^2=k^2=ijk=-1\},\]
with multiplication
\[qq'=(aa'-bb'-cc'-dd')+(ab'+a'b+cd'-dc')i+(ac'+a'c+db'-bd')j+(ad'+a'd+bc'-cb')k.\]
The real and imaginary parts are $\text{Re}q=a, \text{Im}q=bi+cj+dk$. The conjugate and norm are
\[\bar{q}=a-bi-cj-dk,~|q|^2=q\bar{q}=a^2+b^2+c^2+d^2\] satisfying \[\overline{qq'}=\overline{q'}\bar{q},~ |qq'|=|q||q'|.\]
$\bb{H}\simeq \bb{C}^2$ as $q=(a+bi)+(c+di)j$, actually, it can be represented by 2-order complex matrix \[A_q=\left(
                                                            \begin{array}{cc}
                                                              a+bi    & c+di  \\
                                                              -(c-di) & a-bi \\
                                                            \end{array}
                                                          \right),\]
which preserves multiplications, i.e. $A_{qq'}= A_qA_{q'}$. So, if we use notation
$q=(\zeta^1,\zeta^2) (\zeta^1,\zeta^2\in \bb{C})$ for $q=\zeta^1+\zeta^2j$ and take $q'=(\eta^1,\eta^2)$, then
$\bar{q}=(\overline{\zeta^1}, -\zeta^2)$ and
\begin{align*}
qq'&=(\zeta^1\eta^1-\zeta^2\overline{\eta^2},\zeta^1\eta^2+\zeta^2\overline{\eta^1}),\nonumber\\ q\overline{q'}&=(\zeta^1\overline{\eta^1}+\zeta^2\overline{\eta^2},\zeta^2\eta^1-\zeta^1\eta^2).\end{align*}
 For $q=(q_1,\ldots,q_n)\in \bb{H}^n$, we use scalar product $
 q\cdot \overline{q'}=\sum_{j=1}^n q_j\overline{q'_j}$, then quaternionic unitary group, also called compact symplectic group,
 \[Sp(n)\triangleq \{A\in M_n(\bb{H}):~\overline{A^T}A=A\overline{A^T}=I\}\] preserves this product by right action on $\bb{H}^n$.

Quaternionic Heisenberg group is a 2-step nilpotent Lie group identified with its Lie algebra $G=\bb{H}^n\times \text{Im}\bb{H}$ with group multiplication law
\[ uu'=(q,w)(q',w')=(q+q',w+w'+2\text{Im}q\cdot \overline{q'}),\] where $u=(q,w), u'=(q',w')$ are two group elements.
Here $``\cdot"$ is the scalar product in $\bb{H}^n$. This group can be viewed as the nilponent part of Iwasawa
decomposition of Lorentz group $Sp(n+1,1)$, which is the isometry group of the rank one hyperbolic symmetric
space over quaternions, which can be identified with homogeneous space $Sp(n+1,1)/Sp(n+1)\times Sp(1)$. Note here that we need take a little care to choose a proper inner product on Lie algebra to make $G$ defined above to be a group of H-type \cite{k}.
We denote the homogenous norm on $G$ \[ |u|=|(q,w)|=(|q|^4+|w|^2)^{\frac{1}{4}}, \] then the group distance of
two elements is defined naturally by
\[d_G(u,v)= |v^{-1}u| = \left(|q-q'|^4+|w-w'+2\im q\cdot\overline{q'}|^2\right)^{\frac{1}{4}}\] with $u=(q,w),v=(q',w')$.
We use notation $Q=4n+6$ for the homogeneous dimension, $\delta u=(\delta q, \delta^2w)$ for group dilation
and $du=dqdw$ (Lebesgue measure) for Haar measure. The left-invariant vector fields corresponding to one-parameter coordinates subgroups are given by
\begin{align*}
X_j^0=&\frac{\partial}{\partial q_j^0}+2q_j^1\frac{\partial}{\partial w_1}+2q_j^2\frac{\partial}{\partial w_2}+2q_j^3\frac{\partial}{\partial w_3},\\
X_j^1=&\frac{\partial}{\partial q_j^1}-2q_j^0\frac{\partial}{\partial w_1}-2q_j^3\frac{\partial}{\partial w_2}+2q_j^2\frac{\partial}{\partial w_3},\\
X_j^2=&\frac{\partial}{\partial q_j^2}+2q_j^3\frac{\partial}{\partial w_1}-2q_j^0\frac{\partial}{\partial w_2}-2q_j^1\frac{\partial}{\partial w_3},\\
X_j^3=&\frac{\partial}{\partial q_j^3}-2q_j^2\frac{\partial}{\partial w_1}+2q_j^1\frac{\partial}{\partial w_2}-2q_j^0\frac{\partial}{\partial w_3},\\
T^k=&\frac{\partial}{\partial w_k},\qquad \qquad 1\le j\le n,~1\le k\le 3,
\end{align*}
all of which form a basis of Lie algebra of quaternionic Heisenberg group $G$ and the second-order left-invariant differential operator,
sublaplacian, is defined by \[\ca{L}=-\frac{1}{4}\sum_{1\le j\le n, 0\le k\le 3} (X_j^k)^2,\] which is independent of the choice of basis and hypoelliptic from famous theorem of H\"{o}mander.
The fundamental solution of $\ca{L}$ was proved to be
\begin{align*}
\ca{L}^{-1}(\zeta,\eta)= & \frac{2^{\frac{Q}{2}-5}\Gamma(\frac{Q-2}{4})\Gamma(\frac{Q-6}{4})}
{\pi^{\frac{Q}{2}-1}}
d^{2-Q}_G(\zeta,\eta)\\
=& \frac{2^{2n-2}\Gamma(n+1)\Gamma(n)}{\pi^{2n+2}}d^{2-Q}_G(\zeta,\eta).
\end{align*}
An analogue was first proved in \cite{f} for Heisenberg group, and then proved  for
general H-type in \cite{k} when Kaplan introduce the new concept ``H-type".

\emph{Cayley transform and Quaternionic Sphere:}
we denote the quaternionic sphere 
\[S=\{\zeta=
(\zeta',\zeta_{n+1})\in \bb{H}^n\times \bb{H}: |\zeta'|^2+|\zeta_{n+1}|^2=1\}\cong\bb{S}^{4n+3}\]
endowed with Lebesgue sphere measure $d\zeta$. The quaternionic Heisenberg group $G$ is then equivalent to the punctured sphere (quaternionic sphere $S$ minus south pole $o=(0,\ldots,0,-1)$), through ``boundary" Cayley transform defined by
\begin{align}
\ca{C}:~~\qquad G ~&\longrightarrow~ S\setminus \{o\}\nonumber\\
u=(q,w) ~&\longmapsto~ \zeta=(\zeta',\zeta_{n+1})=\left(\frac{2q}{1+|q|^2-w},\frac{1-|q|^2+w}{1+|q|^2-w}\right)\nonumber\\
\ca{C}^{-1}: ~S\setminus \{o\} ~&\longrightarrow~ G\nonumber\\
\zeta=(\zeta',\zeta_{n+1}) ~&\longmapsto~ u=(q,w)=\left(\frac{\zeta'}{1+\zeta_{n+1}}, -\text{Im}\frac{1-\zeta_{n+1}}{1+\zeta_{n+1}}\right)
\end{align} with Jacobian determinant
\begin{align*} |J_{\ca{C}}(u)|= & 2^{Q-3}\big((1+|q|^2)^2+|w|^2\big)^{-\frac{Q}{2}}\\
= & 2^{-3}|1+\zeta_{n+1}|^{Q}.\end{align*}
Here all quotients we use are the left quotient $\frac{q}{q'}=q'^{-1}\cdot q$.
This ``boundary" Cayley transform is a generalization of stereographic projection in Euclidean space and is conformal (Cayley transform on related symmetric spaces is 1-quasiconformal). For Cayley transform on H-type group, see like \cite{iv}.
We define the sphere distance  on $S$ from \cite{ban} in a unifying H-type setting,
\[ d_{S}(\zeta,\eta)=2^{-\frac{1}{2}}|P_\eta \zeta-\eta|^{\frac{1}{2}} \]
with $P_\eta$ being the orthogonal projection operator onto $T^{(2)}_\eta\oplus \bb{R}\eta$, where $T^{(2)}_\eta$ is the
second part of the orthogonal decomposition of tangent space at point $\eta$.
General $d_S$ is checked to be a distance on the sphere in \cite{ban} for I-type group, equivalently H-type
with so-called $J^2$-conditon. For basic results about H-type group and its relationship with I-type group, see like \cite{k,cdkr,yz}.
Using directly the general formula in \cite{ban}, we get by quaternion coordinates
\beq\label{d} d_S(\zeta,\eta)=2^{-\frac{1}{2}}|1-\zeta\cdot\overline{\eta}|^{\frac{1}{2}}.\eeq
It's interesting to note that there exists the following relation between distances on quaternionic sphere $S$ and quaternionic Heisenberg group $G$,
\begin{align}\label{dis}
d_S(\zeta,\eta)=& ((1+|q|^2)^2+|w|^2)^{-\frac{1}{4}}((1+|q'|^2)^2+|w'|^2)^{-\frac{1}{4}}d_G(u,v)\nonumber\\
=& 2^{\frac{3}{Q}-1}|J_{\ca{C}}(u)|^{\frac{1}{2Q}}|J_{\ca{C}}(v)|^{\frac{1}{2Q}}d_G(u,v),
\end{align}
which can be checked
directly, see also \cite{acdb} for general I-type group.
If we take complex coordinates $\zeta'=(\zeta^1_1,\zeta^2_1;\cdots;\zeta^1_n,\zeta^2_n),
\zeta_{n+1}=(\zeta^1_{n+1},\zeta^2_{n+1}), \zeta=(\zeta^1,\zeta^2)$ with
$\zeta^i=(\zeta^i_1,\ldots,\zeta^i_{n+1}), \zeta^i_j\in\bb{C},  1\le i\le 2, 1\le j\le n+1$,
then from (\ref{d})
\beq\label{d-} d_S(\zeta,\eta)=2^{-\frac{1}{2}}(|1-\zeta\cdot_\bb{C}\overline{\eta}|^2+|\zeta\cdot_\bb{C}\sigma(\eta)|^2)^{\frac{1}{4}},\eeq
where $\zeta\cdot_\bb{C}\eta=\sum_{1\le i\le 2,1\le j\le n+1} \zeta_j^i\eta_j^i$ is the complex product and $\sigma(\eta)=(-\eta^2,\eta^1)$ with $ \eta=(\eta^1,\eta^2)$. 
Note that, from (\ref{d-}), we see the distance is not invariant under the action of complex unitary group $U(2n+2)$, so the
associated integral operator with kernel $d_S(\zeta,\eta)$ is not diagonal with respect to the complex bigraded spherical harmonics
irreducible decomposition of $L^2(S)$. However, this is not strange as the quaternion unitary group $Sp(n+1)$ is much
smaller than complex unitary group $U(2n+2)$. Actually, the maximal compact subgroup (stabilize the origin point) of
the isometry group $Sp(n+1,1)$ of the quaternionic (ball-model) symmetric space is $Sp(n+1)\times Sp(1)$,
which leaves invariant and acts transitively on the boundary sphere and
the subsubgroup leaving the north pole fixed is isomorphic to $Sp(n)\times Sp(1)$, then we can realize the spherical principle series representation of
$Sp(n+1)\times Sp(1)$ on $L^2(S)\cong L^2(Sp(n+1)\times Sp(1)/Sp(n)\times Sp(1))$. We need to decompose further real spherical
harmonic irreducible subspaces. 
In other words, we get the Peter-Weyl decomposition of $L^2(S)$ under the representation action of $Sp(n+1)\times Sp(1)$.
Using the explict Jacobi polynomial formula (or hypergeometric series form) for reproducing kernel, we will give the quaternionic
analogue of Funk-Hecke formula, which is useful in proving a critical bilinear inequality for sharp HLS inequality by computing the eigenvalues of corresponding quadratic forms.

The sublaplacian and ``conformal sublaplacian" on $S$ are defined by
\[\ca{L}'=-\frac{1}{4}\sum_{1\le j\le n, 0\le k\le 3} (|J_\ca{C}|^{-\frac{1}{Q}}\ca{C}_*X_j^k)^2,\]
where $\ca{C}_*$ is the induced tangent map of Cayley transform and
\[\ca{D}=\ca{L}'+n(n+1),\] which is an analogue of Geller-type sublaplacian \cite{geller1980laplacian} on Heisenberg group.
The fundamental solution of $\ca{D}$ is given by
\begin{align*}
\ca{D}^{-1}(\zeta,\eta)= & \frac{\Gamma(\frac{Q-2}{4})\Gamma(\frac{Q-6}{4})}
{2^{\frac{Q}{2}}\pi^{\frac{Q}{2}-1}}
d^{2-Q}_S(\zeta,\eta)\\
=& \frac{\Gamma(n+1)\Gamma(n)}{2^{2n+3}\pi^{2n+2}}d^{2-Q}_S(\zeta,\eta)
\end{align*}
and from the definition, we have relation between  sublaplacian $\ca{L}$ on $G$ and conformal sublaplacian $\ca{D}$ on $S$,
\[\ca{L}\left((2|J_{\ca{C}}|)^{\frac{Q-2}{2Q}}(F\circ\ca{C})\right)
=(2|J_{\ca{C}}|)^{\frac{Q+2}{2Q}}(\ca{D}F)\circ\ca{C},\] for all $F\in C^\infty(S)$, see for example \cite{acdb}.
The sublaplacians have one kind of generalization: intertwining operators.

\emph{Intertwining operators on $G$ and $S$:}
general intertwining operators are defined for principle (complementary) series representations of semisimple groups (of real rank one),
here we concern the Lorentz group $Sp(n+1,1)$. See \cite{ks,jw,cow}. That of $SU(n+1,1)$ was also studied by \cite{bfm}
for Heisenberg groups in a geometry and analysis language.

Denote $Aut(G)$ the set of all conformal
transformations on $G$ (diffeomorphisms preserving the contact structure, also called quaternionic automorphisms, that's why we use the notation), element of which is composition of translations, rotations (on $q$),
dilations and inversion
\beq\label{inversion}\sigma_{inv}:(q,w)\mapsto (-\frac{q}{|q|^2-w},-\frac{w}{|q|^4+|w|^2}),\eeq
extension of which on related rank one symmetric space (Siegel domain model) is an isometry. It's a celebrated theorem that $Aut(G)\backsimeq Sp(n+1,1)$ and for inversion on H-type or I-type group, see \cite{kor,ck,cdkr}.
For quaternionic sphere, we denote similarly the set of all conformal transformations on $S$ by
\[Aut(S)\triangleq\{\tau=\ca{C}\circ\sigma\circ\ca{C}^{-1}:\sigma\in Aut(G)\}.\]  Before going on, it's necessary to define Folland-Stein-Sobolev space of exponent $(d,2) (d>0)$, denoted by $W^{d,2}$ here, to be the completion of $C^\infty(S)$ (or $\scr{D}(G)$) w.r.t the norm $\|f\|_{W^{d,2}}=\|(I+\ca{L'})^{\frac{d}{2}}f\|_2$ (or $\|(I+\ca{L})^{\frac{d}{2}}f\|_2$). $\scr{D}(G)$ means the smooth function space of compact support on $G$.

Then for $d\in(0,Q)$, we define intertwining operators $\ca{A}_d$ on $S$ to be any operator satisfying
\[|J_\tau|^{\frac{Q+d}{2Q}}(\ca{A}_dF)\circ\tau=\ca{A}_d\left(|J_\tau|^{\frac{Q-d}{2Q}}(F\circ\tau)\right)\]
for all $F\in C^\infty(S), \tau\in Aut(S)$ and $|J_\tau|$ is the Jacobian determinant of $\tau$. The definition states that operator $\ca{A}_d$
intertwines with two principle series representations $\pi_d,\pi_{-d}$ of $Sp(n+1,1)$, with $\pi_d(\tau): F\rightarrow |J_\tau|^{\frac{Q+d}{2Q}}F\circ \tau$.
It was first proved in \cite{kostant1969existence,jw} (with representation language) that the intertwining operator is diagonal w.r.t to the $Sp(n+1)\times Sp(1)-$irreducible bi-spherical harmonic decomposition
$L^2(S)=\bigoplus_{j\ge k\ge 0} V_{j,k}$ (see (\ref{dec}) in subsection \ref{fh})
and its spectrum is, modula a constant dependent of $d$, uniquely given by,
\beq\label{Ad}
\lambda_{j,k}(\ca{A}_d|_{V_{j,k}})= \frac{\Gamma(j+\frac{Q+d}{4})}{\Gamma(j+\frac{Q-d}{4})}\frac{\Gamma(k+\frac{Q+d}{4}-1)}{\Gamma(k+\frac{Q-d}{4}-1)}.
\eeq
We can also duplicate the calculus in Appendix A of \cite{bfm} to get this spectral decomposition.
Then the operator can be extended onto the Folland-Stein-Sobolev spaces $W^{\frac{d}{2},2}(S)$.
If we choose the $d-$dependent constant to be 1, then by (1) in Theorem \ref{t-eig}~(corresponds to $\alpha=\frac{Q-d}{4}$), the fundamental solution of $\ca{A}_d$ is given by
\beq\label{ad} \ca{A}_d^{-1}(\zeta,\eta)=c_d'd_S^{d-Q}(\zeta,\eta)\eeq
with
\beq\label{c'd}c_d'^{-1}=\frac{2^{\frac{Q-d}{2}+1}\pi^{\frac{Q}{2}-1}\Gamma(\frac{d}{2})}
{\Gamma(\frac{Q-d}{4})\Gamma(\frac{Q-d}{4}-1)},\eeq
and this result is also listed in \cite{acdb} and implicitly contained in \cite{jw}.
The corresponding intertwining operators $\ca{L}_d$ on $G$ can be defined similarly to be any operator satisfying
\[|J_\sigma|^{\frac{Q+d}{2Q}}(\ca{L}_df)\circ\sigma=\ca{L}_d\left(|J_\sigma|^{\frac{Q-d}{2Q}}(f\circ\sigma)\right)\]
for all $\sigma\in Aut(G), f\in \scr{D}(G)$, or equivalently by Cayley transform
\[\ca{L}_d\left((2|J_{\ca{C}}|)^{\frac{Q-d}{2Q}}(F\circ\ca{C})\right)
=(2|J_{\ca{C}}|)^{\frac{Q+d}{2Q}}(\ca{A}_dF)\circ\ca{C},\] for all $F\in\scr{D}(S)$.
Take $Sp(n)$-spherical functions $\{\Phi_k^\lambda\}_{k\in\bb{N},\lambda\in \im\bb{H}\cong \bb{R}^3}$ to be normalized joint radial eigenfunctions of sublaplacian $\ca{L}$ and $\frac{\partial}{\partial w_1}, \frac{\partial}{\partial w_2}, \frac{\partial}{\partial w_3}$, with joint spectrum $\big(2(k+n)|\lambda|, \lambda\big)$ and their formulas are given by
\[\Phi^\lambda_k(q,w)=e^{i\lambda\cdot_{\bb{R}} w-|\lambda||q|^2} L_k^{2n-1}(|\lambda||q|^2),\] where $\lambda\cdot_\bb{R} w=\sum_{j=1}^3\lambda_j w_j$ and $L_k^{2n-1}$ is the $k$-th classical Laguerre polynomial of order $2n-1$. Here, for ``radial", we mean the function depends only on $|q|$ and $w$.
Then we can define ``group Fourier transform" (usually called spherical Fourier transform) on $G$, which acts on radial function in $L^1(G)$ as
\[\hat{}: f \longmapsto \hat{f}(\lambda,k)=\int_G f(u)\Phi_k^\lambda(u).\]
From the intertwining relation and using Fourier transform, we have that the operator is, modula constant multiple, uniquely given by
\beq\label{Ld}
\widehat{\ca{L}_df}(\lambda,k)= |2\lambda|^{\frac{d}{2}}\frac{\Gamma(k+\frac{Q+d}{4}-1)}{\Gamma(k+\frac{Q-d}{4}-1)} \hat{f}(\lambda,k),
\eeq
and ignoring the $d-$dependent constant, it can be explicitly expressed by spectral calculus
\[\ca{L}_d=|2\ca{T}|^{\frac{d}{2}}\frac{\Gamma(|2\ca{T}|^{-1}\ca{L}+\frac{2+d}{4})}{\Gamma(|2\ca{T}|^{-1}\ca{L}+\frac{2-d}{4})}.\]
Then the operator can also be extended to the Folland-Stein-Sobolev spaces $W^{\frac{d}{2},2}(G)$ as before.
The fundamental solution of $\ca{L}_d$ is given by \beq\label{ld}\ca{L}_d^{-1}(u,v)=c_d d_G^{d-Q}(u,v)\eeq with \beq\label{cd}c_d=2^{Q-d-3}c'_d,\eeq where $c_d'$  is defined in (\ref{c'd}).
The operator for related rank one semisimple Lie group was studied in \cite{cow}.

The intertwining operators is one kind of generalization of fractional (conformal) (sub-)laplacian as $\ca{A}_2=\ca{D}, \ca{L}_2=\ca{L}$, and especially in the Euclidean space and sphere, $\ca{L}_d=(-\Delta)^{\frac{d}{2}}$ and $\ca{A}_d$ is the spherical picture obtained from $(-\Delta)^{\frac{d}{2}}$ by stereographic projection. We can also try to consider intertwining operator at $d=Q$, then we may get corresponding analogue of Beckner-Onofri inequality, which is not included in this paper.
Above results were all discussed for Heisenberg group (and complex sphere) in \cite{bfm} and can also be derived from the theory of Knapp-Stein intertwining operators both in the compact and noncompact pictures in \cite{ks}.

\section{Main Results}\label{mr}
Return to HLS inequality, we now can state our main sharp HLS and related HLS-type theorems. In this paper, we use notation $``\sim"$ for equality, modulo a constant multiple.
\subsection{Sharp HLS Inequalities}\label{mr1}
\bthm\label{t-q}\emph{[{Sharp HLS on Quaternionic Heisenberg Group}]}\\
Let $4\le \lambda<Q=4n+6, p=\frac{2Q}{2Q-\lambda}$, then $\forall ~f,g\in L^p(G)$,
\beq\label{hls}
\Big|\iint_{G\times G}\frac{\overline{f(u)}g(v)}{|u^{-1}v|^\lambda}dudv\Big| \le C_\lambda \|f\|_p\|g\|_p
\eeq with sharp constant
\begin{align}\label{C}
C_\lambda
=&
\frac{2^{1-\frac{2n\lambda}{Q}}\pi^{2n+2}\Gamma(\frac{Q-\lambda}{2})|S|^{1-\frac{2}{p}}}{\Gamma(\frac{2Q-\lambda}{4})\Gamma(\frac{2Q-\lambda}{4}-1)}\nonumber\\
=& \left(\frac{\pi^{\frac{Q-2}{2}}}{2^{\frac{Q-8}{2}}}\right)^{\frac{\lambda}{Q}}
\frac{((\frac{Q-4}{2})!)^{1-\frac{\lambda}{Q}}\Gamma(\frac{Q-\lambda}{2})}{\Gamma(\frac{2Q-\lambda}{4})\Gamma(\frac{2Q-\lambda}{4}-1)},
\end{align} where $|S|=\frac{2\pi^{2n+2}}{(2n+1)!}$ is the surface area of $S$.
Moreover, all extremizers are given by
\beq f\sim g \sim (|J_{\ca{C}}\circ\sigma||J_\sigma|)^{\frac{1}{p}} 
\sim\left||q|^2+w-2q_0\cdot \bar{q}+r_0\right|^{-\frac{2Q-\lambda}{2}},\eeq
with $\sigma\in Aut(G), 
q_0\in \bb{H}^n, r_0\in \bb{H}$, satisfying $\emph{Re}r_0>|q_0|^2$, and we can choose one
$\sigma=\ca{S}_{\delta_0}\circ \ca{L}_{u_0}$ or ~$\ca{S}_{\delta}\circ\ca{C}^{-1}\circ A_\xi\circ \ca{C}$, where $\ca{S}_{\delta_0},\ca{S}_\delta$ are dilations, $\ca{L}_{u_0}$ is left translation and
$A_\xi$ is a rotation in $Sp(n+1)$ s.t. $A_\xi^{-1}(0,\ldots,0,1)=\frac{\xi}{|\xi|}$, with parameters
$\delta_0=(\re r_0-|q_0|^2)^{-\frac{1}{2}}, u_0=(q_0,-\im r_0)$, $\xi=(\frac{2q_0}{r_0+1},\frac{r_0-1}{r_0+1}), \delta=\sqrt{\frac{1\mp|\xi|}{1\pm|\xi|}}$.
\ethm

Using the Cayley transform, we can give the sphere editon of last theorem: with the relation (\ref{dis}) between two distances on $G$ and $S$ , we have the following equivalent sphere editon of sharp HLS inequality above,
through the correspondence between functions $f$ on $G$ and $\tilde{f}$ on $S$,
\beq\label{corr}\tilde{f}(\zeta)=f(\ca{C}^{-1}\zeta)|J_{\ca{C}^{-1}}|^{\frac{1}{p}}.\eeq
\bthm\label{t-s} \emph{[{Sharp HLS on Quaternionic Sphere}]}\\
Let $4\le \lambda<Q=4n+6, p=\frac{2Q}{2Q-\lambda}$, then $\forall ~f,g\in L^p(S)$,
\beq\label{shls}
\Big|\iint_{S\times S}\frac{\overline{f(\zeta)}g(\eta)}{d_S^\lambda(\zeta,\eta)}d\zeta d\eta\Big| \le C'_\lambda \|f\|_p\|g\|_p
\eeq with sharp constant
\begin{align}\label{C'}
C'_\lambda=&2^\frac{(4n+3)\lambda}{Q}C_\lambda\nonumber\\
=&
\frac{2^{1+\frac{\lambda}{2}}\pi^{2n+2}\Gamma(\frac{Q-\lambda}{2})|S|^{1-\frac{2}{p}}}{\Gamma(\frac{2Q-\lambda}{4})\Gamma(\frac{2Q-\lambda}{4}-1)}\nonumber\\
=& \left(\frac{2\pi^{\frac{Q-2}{2}}}{(\frac{Q-4}{2})!}\right)^{\frac{\lambda}{Q}}
\frac{2^{\frac{\lambda}{2}}(\frac{Q-4}{2})!\Gamma(\frac{Q-\lambda}{2})}{\Gamma(\frac{2Q-\lambda}{4})\Gamma(\frac{2Q-\lambda}{4}-1)},
\end{align} where $|S|=\frac{2\pi^{2n+2}}{(2n+1)!}$ is the surface area of $S$. Moreover, all extremizers are given by
\beq f\sim g\sim |J_\tau|^{\frac{1}{p}}\sim|1-\xi\cdot\bar{\zeta}|^{-\frac{2Q-\lambda}{2}},\eeq
with $\tau\in Aut(S), \xi\in \bb{H}^{n+1}, |\xi|<1$, and we can choose one $\tau=\ca{C}\circ \ca{S}_{\delta_0}\circ \ca{L}_{u_0}\circ \ca{C}^{-1}$ or ~$\ca{C}\circ\ca{S}_\delta\circ\ca{C}^{-1}\circ A_\xi$,
where $\ca{S}_{\delta_0},\ca{S}_\delta$ are group dilations, $\ca{L}_{u_0}$ is group left translation, with parameters $\delta_0=\frac{|1+\xi_{n+1}|}{\sqrt{1-|\xi|^2}}, u_0=(\frac{\xi'}{1+\xi_{n+1}},-\im \frac{1-\xi_{n+1}}{1+\xi_{n+1}}), \delta=\sqrt{\frac{1\pm|\xi|}{1\mp|\xi|}}$, and $A_\xi$ is a rotation in $Sp(n+1)$ s.t.
$A_\xi^{-1}(0,\ldots,0,1)=\frac{\xi}{|\xi|}$.
\ethm
At first, we give several small remarks about the two theorems.

\textbf{Remark:}
\bitm
\item
Existence of extremizers holds for all $0<\lambda<Q$ for Theorem \ref{t-q} and \ref{t-s}.
Several a bit standard methods can be used to prove the existence while compactness is the basic idea.
We will still borrow the main argument from \cite{fl} and give the outline proof.

\item
The large conformal symmetry group of the HLS inequality (\ref{hls}) consists of constant multiples, left-translations, dilations,
rotations (on $q$ variable) and inversion: $f(u)\mapsto f(\sigma_{inv} u) |u|^{-\frac{2Q}{p}}$, where $\sigma_{inv}$ is the
inversion transform (\ref{inversion}). In other words, the inequality is invariant under the conformal action
$f\mapsto f\circ\sigma |J_\sigma|^{\frac{1}{p}}, \forall \sigma\in Aut(G)$. Similarly, HLS inequality on sphere (\ref{shls}) is invariant under conformal action $f\mapsto f\circ\tau |J_\tau|^{\frac{1}{p}}, \forall \tau\in Aut(S)$.
\item
In Theorem \ref{t-q}, modulo constant multiples, left-translations and dilations,
extremizer exists uniquely,  i.e. sharp equality (\ref{hls}) holds if and only if
\[ f=g=\big((1+|q|^2)^2+|w|^2\big)^{-\frac{2Q-\lambda}{4}}.\]
\item
In Theorem \ref{t-s}, modulo constant multiples, quaternionic rotations, group dilations and
Cayley transform, extremizer exists uniquely, i.e. sharp equality (\ref{shls}) holds if and only if
\[ f=g=1.\]
\item
For finding the explicit formula for extremizers, we select a specially chosen extremizer which is proved to be constant for
$\lambda> 4$ and be in direct sum of subspaces $V_{0,0}\bigoplus_{j\ge k\ge 2}V_{j,k}$ (see (\ref{dec}))
for $\lambda=4$.
Fortunately, for $\lambda=4$, we can recur to the Euler-Lagrange equation of first variation
to restrict further the extremizers to be in $V_{0,0}$ which consists of constant functions.
We can also see from the original functional of inequality. For $\lambda<4$, we have no that good result, but prove weakly that constant is ``local" extremizer, i.e., there exists a small domain around constant function w.r.t the Lebesgue norm, s.t., the inequality functional take maximum in the center.
\eitm
We leave the proofs of the two main theorems in section \ref{p-hls} and want to give some related sharp inequalities first.
\subsection{Sharp Sobolev-type Inequalities}\label{mr2}
We can also write above sharp HLS inequalities in a dual form concerning intertwining operators, which we may call sharp
conformally-invariant Sobolev-type inequalities, noticing that the fundamental solution of intertwining operators $\ca{A}_d$ on $S$ ($\ca{L}_d$ on $G$)  are constant multiple of $d_S(\zeta,\eta)^{d-Q}(d_G(u,v)^{d-Q})$, see (\ref{ad}) and (\ref{ld}), here $d\in(0,Q)$.
\bthm\label{cis} \emph{[Sharp Conformally-Invariant Sobolev Inequality]}\\
\emph{(1)} Let $0<d\le Q-4, p=\frac{2Q}{Q-d}$, then $\forall f\in W^{\frac{d}{2},2}(G)$,
\beq
\int_G \bar{f}\ca{L}_df\ge \tilde{C}_d\|f\|_p^2,
\eeq with sharp constant
\[\tilde{C}_d=(c_dC_{Q-d})^{-1},\]
 and all extremizers
\[f\sim \left||q|^2+w-2q_0\cdot \bar{q}+r_0\right|^{-\frac{Q-d}{2}},\] with $q_0\in \bb{H}^n, r_0\in \bb{H}$,
satisfying $\emph{Re}r_0>|q_0|^2$.\\
\emph{(2)} Let $0<d\le Q-4, p=\frac{2Q}{Q-d}$, then $\forall f\in W^{\frac{d}{2},2}(S)$,
\beq
\int_S \bar{f}\ca{A}_df\ge \tilde{C_d'} \|f\|_p^2,
\eeq with sharp constant
\[\tilde{C_d'}=(c'_dC'_{Q-d})^{-1},\]
 and all extremizers
\[f\sim |1-\xi\cdot\bar{\zeta}|^{-\frac{Q-d}{2}},\] with $\xi\in \bb{H}^{n+1}, |\xi|<1$.\\
Constant $c_d, c'_d$ are defined by \emph{(\ref{cd})(\ref{c'd})} in preliminary section, and $C_{Q-d},C'_{Q-d}$ is given by \emph{(\ref{C}),(\ref{C'})} in Theorem \emph{\ref{t-q},\ref{t-s}}.
\ethm

It's obvious that we only need to consider the inequality for real-valued, or even nonnegative functions. Besides, we can easily checked that the inequalities are invariant under the action of conformal group, i.e., under the action of transformations $f\mapsto f\circ\sigma |J_\sigma|^{\frac{1}{p}}, \forall \sigma\in Aut(G) ~(\text{or}~ Aut(S))$. Also, as in \cite{fl}, we can give the endpoint limit case of sharp HLS inequality at $\lambda=Q$, using standard functional limit argument.
The endpoint case corresponds to Log-Sobolev inequality. We list it in the sphere framework.

\bthm\label{end} \emph{[Sharp Log-Sobolev Inequality]}\\
$\forall f\ge 0\in L^2\emph{Log}L(S)$, normalized by $\int_S f^2=|S|$,
\beq\label{log}\iint_{S\times S}\frac{|f(\zeta)-f(\eta)|^2}{d^Q_S(\zeta,\eta)}d\zeta d\eta \ge C \int_S f^2\log f^2
, \eeq
with sharp constant \[C=
\frac{2^{\frac{Q}{2}+3}\pi^{\frac{Q}{2}-1}}{Q\Gamma(\frac{Q}{4}-1)\Gamma(\frac{Q}{4})}
,\] and some extremizers
\[ f\sim |1-\xi\cdot\bar{\zeta}|^{-\frac{Q}{2}},\] satisfying normalized condition with nonzero $\xi\in\bb{H}^{n+1}, |\xi|<1$.
\ethm
Now we give a simple proof of Theorem \ref{cis} and \ref{end}.
First, we assume the validity of Theorem \ref{t-q} and \ref{t-s}.\\
\emph{Proof of Theorem} \ref{cis}:
\bpf
Because of relationship between $\ca{L}_d$ and $\ca{A}_d$ and similar arguments, we only prove on group for simplicity.
By Plancherel formula, Cauchy-Schwartz inequality and properties (\ref{Ld}) and (\ref{ld}) of intertwining operators, we have
\begin{align*}
\ip{f,g}^2=\ip{\hat{f},\hat{g}}^2=&\ip{\left(|2\lambda|^{\frac{d}{2}}\frac{\Gamma(k+\frac{Q+d}{4}-1)}{\Gamma(k+\frac{Q-d}{4}-1)}\right)^{\frac{1}{2}}\hat{f},\left(|2\lambda|^{\frac{d}{2}}\frac{\Gamma(k+\frac{Q+d}{4}-1)}{\Gamma(k+\frac{Q-d}{4}-1)}\right)^{-\frac{1}{2}}\hat{g}}^2\\
\le&\ip{\hat{f},\widehat{\ca{L}_df}}\ip{\hat{g},\widehat{\ca{L}^{-1}_dg}}\\
=&\ip{\hat{f},\widehat{\ca{L}_df}}\ip{\hat{g},c_d \widehat{|u|^{d-Q}*g}}.
\end{align*}
So we have the following equivalent form
\[\abs{\ip{f,g}}^2\le c_d \ip{f,\ca{L}_df} \ip{g,|u|^{d-Q}*g},\] i.e.
\[\abs{\int_G \bar{f}g}^2 \le c_d \int_G \bar{f}\ca{L}_df \iint_{G\times G} \frac{\overline{g(u)}g(v)}{d_G^{Q-d}(u,v)}dudv. \]
So, from sharp HLS inequality in Theorem \ref{t-q}, we have
\begin{align*}
\|f\|^2_p \le & c_d \int_G \bar{f}\ca{L}_df \sup_{\|g\|_{p'}=1}\iint_{G\times G} \frac{\overline{g(u)}g(v)}{d_G^{Q-d}(u,v)}dudv\\
=&c_d C_{Q-d} \int_G \bar{f}\ca{L}_df,
\end{align*} and ``=" holds if and only if $g$ is extremizer for sharp HLS inequality and $g\sim f^{p-1}$, i.e. $f$ can only be the form in the Theorem.
Actually, from the obvious opposite direction, we can see the equivalence of sharp conformally-invariant Sobolev inequality and HLS inequality.
\epf
\emph{Proof of Theorem} \ref{end}:
\bpf
From the normalized assumption,
\[\iint_{S\times S} \frac{f^2(\zeta)+f^2(\eta)}{d^\lambda_S(\zeta,\eta)}d\zeta d\eta=2C'_\lambda|S|^{\frac{2}{p}}.\]
Subtracting two mutiple of inequality (\ref{shls}) and taking limitation $\lambda\rightarrow Q$, we get
\begin{align*}
\iint_{S\times S}\frac{|f(\zeta)-f(\eta)|^2}{d^Q_S(\zeta,\eta)}d\zeta d\eta
\ge& \lim_{\lambda\rightarrow Q}2C'_\lambda(|S|^{\frac{2}{p}}-\|f\|_p^2)\\
=&\frac{2^{\frac{Q}{2}+3}\pi^{2n+2}}{\Gamma(\frac{Q}{4}-1)\Gamma(\frac{Q}{4})}
\lim_{\lambda\rightarrow Q}\frac{|S|^{\frac{2}{p}}-\|f\|_p^2}{Q-\lambda}.
\end{align*}
It's easy to check
\[
\frac{|S|^{\frac{2}{p}}-\|f\|_p^2}{Q-\lambda}=\frac{p-2}{{Q-\lambda}}\frac{|S|^{\frac{2}{p}}-\|f\|_p^2}{p-2}
\xrightarrow{\lambda \rightarrow Q} \left(-\frac{2}{Q}\right)\left(-\frac{\int_S f^2\log{f^2}}{2}\right)=\frac{\int_S f^2\log{f^2}}{Q}.
\]
Then the theorem is proved after checking equality can be achieved by some extremizers $\sim |1-\xi\cdot\bar{\eta}|^{-\frac{Q}{2}}$,
which is the limit of extremizers for sharp HLS inequality (\ref{shls}). Actually, if we denote $h\sim |1-\xi\cdot\bar{\eta}|^{-\frac{Q}{2}}$,
$h_\lambda\sim |1-\xi\cdot\bar{\eta}|^{-\frac{2Q-\lambda}{2}}$, both satisfying normalized condition, then
$h_\lambda= \left(\frac{|S|}{\int_S h^{\frac{4}{p}}}\right)^{\frac{1}{2}} h^{\frac{2}{p}}, ~h_\lambda\rightrightarrows h (\lambda\rightarrow Q)$ and we have
\[\iint_{S\times S}\frac{|h_\lambda(\zeta)-h_\lambda(\eta)|^2}{d^\lambda_S(\zeta,\eta)}d\zeta d\eta
= 2C'_\lambda(|S|^{\frac{2}{p}}-\|h_\lambda\|_p^2).\]
As $\lambda\rightarrow Q$, the right side of above equality again converges to the right side of (\ref{log}), replacing $f$ by $h$,
while the left side converges to the left side of (\ref{log}), replacing $f$ by $h$ because of uniformly convergence
$h_\lambda\rightrightarrows h$.   Another way to prove sharpness is using counter-example
$f(\zeta)=\sqrt{1-\epsilon} +\epsilon \re\zeta_1$ from \cite{fl}.
\epf

\section{Proof of Sharp HLS Inequalities (Theorem \ref{t-q} and \ref{t-s})}\label{p-hls}
The left space of this paper is given to the proof of our main results of sharp HLS on quaternionic group and sphere. In order to express clearly the idea thread, we split it into 6 steps.\\
\emph{Proof of Theorem} \ref{t-q} $and$ \ref{t-s}:
\subsection{Existence for General $\lambda$}
\emph{Existence of extremizer.}
Arguments for exsitence in \cite{fl} can be easily transformed to our case, considering compactness argument, Fatou lemma and
``refined HLS". Concentration compactness from \cite{lio} is an alternative method, see \cite{h}. We here give an outline of
proof from the method of \cite{fl}. First, we note that the kernel of power of distance is positive, actually, it's the
fundamental solution of intertwining operators $\ca{L}_d^{-1}\sim |u|^{-\lambda}$, with $d=Q-\lambda$.
Take $k\sim \ca{L}_d^{-\frac{1}{2}} \delta_0\sim \ca{L}_d^{-\frac{1}{2}}\ca{L}_{\frac{d}{2}}|u|^{-\frac{Q+\lambda}{2}}$,
with $\delta_0$ the standard Dirac function at zero, then $|u|^{-\lambda}=k\ast k$ and from the $L^p$ boundness of
$\ca{L}_d^{-\frac{1}{2}}\ca{L}_{\frac{d}{2}}$ and Marcinkiewicz interpolation theorem, we get
$k\in L^{\frac{2Q}{Q+\lambda}, \infty}$. $k$ is also even, real-valued, homogeneous of order $-\frac{Q+\lambda}{2}$.
To prove the $L^p$ bound, we use and check the codition for the Marcinkiewicz multiplier theorem in \cite{msr}.
Then the sharp problem changes to that of \beq\label{l2}\|f\ast k\|_q \lesssim \|f\|_2,\quad q=p'=\frac{2Q}{\lambda} ~\text{is the conjugate exponent}.\eeq When using the compactness argument for above inequality, we should be alert as the weak-limit of extremizing sequence is easy to vanish. So, the following ``quaternionic edition" of an enforced estimate called ``refined HLS" inequality kills this bug easily.
\blem\label{ref}
Let $0<\lambda<Q, q=\frac{2Q}{\lambda}, k\sim \ca{L}_d^{-\frac{1}{2}} \delta$, then
\[\|f\ast k\|_q \lesssim \|f\|_2^{\frac{\lambda}{Q}}\left(\sup_{\beta>0}\beta^{\frac{\lambda}{4}}\|e^{-\beta\ca{L}}\beta\ca{L}(f\ast k)\|_\infty\right)^{\frac{Q-\lambda}{Q}},\] with $e^{-\beta\ca{L}}$ the heat semigroup of sublaplacian.
\elem
From lemma \ref{ref} and properties of $k$ and heat kernel, we get a non-zero weak limit of extremizing sequence in $L^2$ for
(\ref{l2}) $f_j\rightharpoonup f \not\equiv 0$, satisfying $f_j\ast k \rightarrow f\ast k, a.e.$ We assume that $\|f_j\|_2=1$,
then from general ``Fatou lemma" (for this useful lemma and application, see \cite{lie} and the references there),
we get $f_j\xrightarrow{L^2} f$,
which tells that $f$ is the extremizer for (\ref{l2}) and therefore (\ref{hls}) and (\ref{shls}).

Lemma \ref{ref} is a corollary of the following lemma (take $r=2, s=\frac{Q-\lambda}{2}$ in the lemma) and the $L^2$ boundness of $\|\ca{L}^{\frac{s}{2}}(f\ast k)\|_2 \lesssim \|f\|_2$, which is again proved from multiplier theorem in \cite{msr}.
\blem
Let $1\le r\le \infty, 0<s<\frac{Q}{r}, q=\frac{rQ}{Q-rs}$, then we have
\[\|f\|_q \lesssim \|f\|^{1-\frac{sr}{Q}}_{\dot{W}^{s,r}(G)}\|f\|_{\dot{B}_{\infty,\infty}^{s-\frac{Q}{r}}}^{\frac{sr}{Q}},\]
with homogeneous Sobolev norm $\|f\|_{\dot{W}^{s,r}(G)}\triangleq \|\ca{L}^{\frac{s}{2}}f\|_r$ and
Besov norm $\|f\|_{\dot{B}_{p,q}^s}$ defined analogously with Euclidean correspondence (see \emph{\cite{fy,sz}}).
\elem
\bpf
We use the heat kernel characterization of homogeneous Besov space in \cite{fy} (there is some error about the condition for order $k$, but we can modify from the proof),
\[\|f\|_{\dot{B}_{\infty,\infty}^{s-\frac{Q}{r}}}\sim \sup_{t>0} t^{-\frac{s-\frac{Q}{r}}{2}}\|(t\ca{L})^k e^{-t\ca{L}}f\|_\infty,\] $\forall k\in \bb{N}^+$ and for ``$\gtrsim$"  we only need $k\in \bb{N}$. 
We write $f$ into the integral of heat flow of sublaplacian operator,
\begin{align*}
f=&\int_0^\infty e^{-t\ca{L}}\ca{L}f dt\\
=&\left(\int_0^A+\int_A^\infty\right)e^{-t\ca{L}}\ca{L}f dt\\
\triangleq & I_1+I_2,
\end{align*} with $A$ a constant to be fixed later.
For the first term, we know that the heat kernel is Schwartz function, furthermore, from the explicit formula (\cite{z} and \cite{yz} for H-type)
 \[p_t(u)\sim \int_{\bb{R}^3} e^{-i<w,w'>_{\bb{R}^3}}e^{-|q|^2|w'|\coth2t|w'|}\left(\frac{|w'|}{\sinh2t|w'|}\right)^{2n}dw',\]
with $u=(q,w)\in G$, we have exponential rapid decreasing estimate for heat kernel
\[|p_t(u)|\sim e^{-O(|q|^2+|w|)},\] when $|u|\rightarrow \infty.$ So, easy argument tells that the kernel $\varphi$ of operator $e^{-t\ca{L}}(-t\ca{L})^{1-\frac{s}{2}}$ satisfies
\[\varphi \in L^1, \quad \psi(u)\triangleq\sup_{|u'|\ge |u|} \varphi(u') \in L^1,\] and using group analogue of classical result about convolution operator (XIII 7.11 in \cite{ste-}), we get estimate
\[\left|e^{-t\ca{L}}(-t\ca{L})^{1-\frac{s}{2}}f\right|\le \|\psi\|_1 Mf, \] where \[M: f\mapsto Mf(u)=\sup_{r>0}\frac{1}{|\{v:|v^{-1}u|<r\}|} \int_{|v^{-1}u|<r} f(v)dv,\] is the maximal function operator and is bounded from $L^r$ to $L^r, \forall r>1$. Then we have
\begin{align*}
|I_1|=&\left|\int_0^A(-t)^{\frac{s}{2}-1}\left(e^{-t\ca{L}}(-t\ca{L})^{1-\frac{s}{2}}\right)\ca{L}^{\frac{s}{2}}f dt\right|\\
\lesssim &A^{\frac{s}{2}} M(\ca{L}^{\frac{s}{2}}f).
\end{align*}
For the second term, from the characterization,
\[|e^{-t\ca{L}}\ca{L}f|\lesssim t^{\frac{s-\frac{Q}{r}}{2}-1}\|f\|_{\dot{B}_{\infty,\infty}^{s-\frac{Q}{r}}},\] we have
\[|I_2|\lesssim A^{\frac{s-\frac{Q}{r}}{2}}
\|f\|_{\dot{B}_{\infty,\infty}^{s-\frac{Q}{r}}}.\]
Combining the estimates for two terms and taking
$A=\left(\frac{\|f\|_{\dot{B}_{\infty,\infty}^{s-\frac{Q}{r}}}}{M(\ca{L}^{\frac{s}{2}}f)}\right)^{\frac{rs}{Q}}$, we have
\[|f|\lesssim {\|f\|^{\frac{rs}{Q}}_{\dot{B}_{\infty,\infty}^{s-\frac{Q}{r}}}}{M^{1-\frac{rs}{Q}}(\ca{L}^{\frac{s}{2}}f)}.\]
Then the lemma is proved using the $L^r-$boundness of maximal operator (\cite{ste-}) and the relation of exponents.
\epf

Now, it suffices to compute the \emph{explicit form of extremizers}. \\
For simply computation, we deal with the sharp inequality on framework of sphere.
From the proof of existence of extremizer, we know that the integral kernal $d_S(\zeta,\eta)^{-\lambda}$ is positive definite.
Therefore we can restrict the sharp problem to $f=g$ case.
By standard argument, we can further restrict the extremizer to be a complex multiple of positive real-value function.
Acutally, for any $f=a+ib$, with $a,b$ respectively the real and image part function on $S$.
Then the left side of the inequality (\ref{shls}) $I(f)=I(a)+I(b)$ because of the symmetry of the distance.
By Cauchy-Schwartz inequality, $I(f)\le I(|f|=\sqrt{a^2+b^2})$ with the equality holds if and only if
$a(\zeta)a(\eta)\equiv b(\zeta)b(\eta)$ and $a(\zeta)a(\eta)\ge 0$, which gives that $f=(\frac{a}{b}+i)b=c|b|, c\in \bb{C}$.
So we can assume the extremizer $h$ is nonnegtive. The vanishing first variation give Euler-Lagrange equation for $h$,
\beq\label{EL} h^{p-1}(\zeta)\sim \int_S \frac{h(\eta)}{|1-\zeta\cdot\bar{\eta}|^{\frac{\lambda}{2}}}d\eta,\eeq
which tells $h$ is positive a.e..
The second variation of functional associated to inequality (\ref{shls}) is
\beq\label{sv} \iint_{S\times S} \frac{\overline{\varphi(\zeta)}\varphi(\eta)}{|1-\zeta \cdot \bar{\eta}|^{\frac{\lambda}{2}}}d\zeta d\eta\int_S h^p-(p-1)\iint_{S\times S} \frac{h(\zeta)h(\eta)}{|1-\zeta\cdot\bar{\eta}|^{\frac{\lambda}{2}}}d\zeta d\eta\int_S h^{p-2}|\varphi|^2\le 0\eeq
for all $\varphi$ satisfying $\int_S h^{p-1}\varphi=0$.

\subsection{Zero Center-Mass Condition}\label{zcm}
The sharp HLS inequality on group is invariant under translation, dilation and of course constant multiple,
while the sphere-edition is preserved by the action of quaternionic unitary group $Sp(n+1)$.
So, through Cayley transform, we can intuitionally try to assume the extremizer satisfies a certain zero center-mass condition,
which breaks the huge conformal symmetry. Indeed, the conformal invariant group plays a critical role here, while often absent for other inequalities.

We define a conformal map on $S\setminus \{A^{-1}o\}$:
\beq \gamma^\delta_A=A^*\circ\ca{C}\circ\ca{S}_\delta\circ\ca{C}^{-1}\circ A\eeq
with any $A\in Sp(n+1),~\ca{S}_\delta, \delta>0$ the dilation on $G$.
Actually, take $\xi=A^{-1}(0,\ldots,0,1)$, then the map $\gamma^\delta_A$ is $(\delta,\xi)$-determinate, but independent of the choice of $A$, so we can use notation $\gamma^\delta_\xi=\gamma^\delta_A$: from diagram \beq\zeta\xrightarrow{A}(\frac{2q}{1+|q|^2-w},\frac{1-|q|^2+w}{1+|q|^2-w})\xrightarrow{\ca{S}_\delta\circ\ca{C}^{-1}} (\delta q,\delta^2 w)\stackrel{\ca{C}}{\rightarrow} (\frac{2\delta q}{1+\delta^2(|q|^2-w)},\frac{1-\delta^2(|q|^2-w)}{1+\delta^2(|q|^2-w)})\xrightarrow{A^*} \gamma^\delta_\xi(\zeta)\eeq
and unitary property of $A$
\[\zeta\cdot \bar{\xi}=A\zeta\cdot(0,\ldots,0,1)=\frac{1-|q|^2+w}{1+|q|^2-w},\]
we give the explicit form
\beq \gamma^\delta_\xi(\zeta)=\frac{2\delta}{1+\zeta\cdot\bar{\xi}+\delta^2(1-\zeta\cdot\bar{\xi})}
 \, \big(\zeta-(\zeta\cdot\bar{\xi})\xi\big)+\frac{1+\zeta\cdot\bar{\xi}-\delta^2(1-\zeta\cdot\bar{\xi})}{1+\zeta\cdot\bar{\xi}+\delta^2(1-\zeta\cdot\bar{\xi})} \, \xi \eeq
with boundary limit
\beq\label{5}\gamma^\delta_\xi(\zeta)\xrightarrow{\delta\rightarrow 0} \xi,\quad \gamma^\delta_\xi(\zeta)\xrightarrow{\delta\rightarrow 1} \zeta.\eeq
Note that as $\delta\rightarrow 1$, the convergence is uniform in $(\xi,\zeta)$ on total $S\times S$ ,
and as $\delta \rightarrow 0$, the convergence is uniform on any slightly small subset
\[E=\{(\zeta,\xi)\in S\times S:~|1+\zeta\cdot\bar{\xi}|\ge \epsilon>0\}.\]

For any $\int_S f=1$, take \[F(r\xi)=\int_S \gamma^{1-r}_\xi(\zeta) f(\zeta)d\zeta,\] defined on unit ball minus origin point, then we have origin and boundary limit uniformly in $\xi$:
\beq F(r\xi)\xrightarrow[\xi]{r\rightarrow 0}
\int_S \zeta f(\zeta)d\zeta, \quad F(r\xi)\xrightarrow[\xi]{r\rightarrow 1} \xi.\eeq
The notation $\xrightarrow[\xi]{r\rightarrow 0}$ means the convergence is uniformly in $\xi$ when $r\rightarrow 0$.
Actually, the $r\rightarrow 0$ case is a direct consequence of the $\delta\rightarrow 1$ case in (\ref{5}),
while the $r\rightarrow 1$ case comes from the $\delta\rightarrow 0$ case in (\ref{5}) and
\begin{align*}
|F(r\xi)-\xi| \le& \big(\int_{E|_S} + \int_{[(S\times S)\setminus E]|_S}\big) |\gamma^{1-r}_\xi(\zeta)-\xi|f(\zeta)d\zeta  \\
        =& (\int_{d_S(\zeta,-\xi)\ge \epsilon}+\int_{d_S(\zeta,-\xi)\le \epsilon})|\gamma^{1-r}_\xi(\zeta)-\xi|f(\zeta)d\zeta\\
        =& I_1+I_2,
\end{align*} with $I_2\xrightarrow[\xi]{\epsilon \rightarrow 0} 0$  and $I_1\xrightarrow[\xi]{r\rightarrow 1} 0$.\\
So $F$ can be extended to be a continous function on closed unit ball in $\bb{H}^{n+1}\simeq \bb{R}^{4n+4}$ and $F(\xi)=\xi$ on
the sphere. By Brouwer's fixed point theorem, there exists at least one zero point, i.e. there exists one $(\delta_0,\xi_0)$ ~s.t.~ $\int_S \gamma^{\delta_0}_{\xi_0}(\zeta)f(\zeta)d\zeta=0$. Take $\gamma=\gamma^{\delta_0}_{\xi_0}$, then we have the following zero center-mass condition result:
\blem\label{lem}
For any positive extremizer of sharp HLS inequality \emph{(\ref{shls})}, there exists a comformal transformation $\gamma$  s.t.
by $h\mapsto \tilde{h}=|J_{\gamma^{-1}}|^{\frac{1}{p}}h\circ \gamma^{-1}$, we get another positive extremizer $\tilde{h}$ satisfying zero center-mass condition:
\beq\label{zero} \int_S \zeta \tilde{h}^p(\zeta)d\zeta =0. \eeq
\elem
\bpf take $f=h^p$ with $h$ the positive extremizer, then from above arguments, there exists one conformal transformation $\gamma$~  s.t.~ $\int_S \gamma(\zeta)h^p(\zeta)d\zeta=0.$
Change the variables, then we get $\int_S \zeta \tilde{h}^p(\zeta)d\zeta=0.$ Besides, from the correspondence (\ref{corr})
between extremizers of sharp HLS on quaternion Heisenberg group and sphere-edition HLS, we know that
$\tilde{h}(\zeta)\sim|J_{\ca{C}}(A\zeta)|^{\frac{1}{p}}|J_{\ca{C}^{-1}}(A\gamma(\zeta))|^{\frac{1}{p}}h\circ\gamma^{-1}(\zeta)$,
which is still an extremizer for (\ref{shls}).
\epf
The zero-center mass condition is critical here in restricting our extremizers into a small, particular function class.

\subsection{Inverse Second Variation Inequality}\label{ss-isv}
We now return to prove the main theorem. Now, from Lemma \ref{lem}, we can assume extremizer $h$
for inequality (\ref{shls}) satisfies (\ref{zero}). Then we will try to prove $h$ can only be constant function.
Substitute $\varphi(\zeta)=h\zeta^1_j, h\zeta^2_j$, both satisfying $\int_S h^{p-1}\varphi=0$,
into the second variation (\ref{sv}) and summing the results,
we get
\beq \iint_{S\times S}\frac{h(\zeta)\sum_{j=1}^{n+1}(\overline{\zeta^1_j}\eta^1_j+\overline{\zeta^2_j}\eta^2_j)h(\eta)}{|1-\zeta\cdot\bar{\eta}|^{\frac{\lambda}{2}}}d\zeta d\eta\le
(p-1)\iint_{S\times S} \frac{h(\zeta)h(\eta)}{|1-\zeta\cdot\bar{\eta}|^{\frac{\lambda}{2}}}d\zeta d\eta.
\eeq
From the symmetry of left side integrand on $(\zeta,\eta)$ and
$\bar{\zeta}\cdot_\bb{C}\eta+\bar{\eta}\cdot_\bb{C}\zeta=\bar{\zeta}\cdot\eta+\bar{\eta}\cdot\zeta$,
we have
\beq\label{sv-}
\iint_{S\times S}\frac{h(\zeta)(\bar{\zeta}\cdot\eta+\bar{\eta}\cdot\zeta)h(\eta)}{|1-\zeta\cdot\bar{\eta}|^{\frac{\lambda}{2}}}d\zeta d\eta\le 2(p-1)\iint_{S\times S} \frac{h(\zeta)h(\eta)}{|1-\zeta\cdot\bar{\eta}|^{\frac{\lambda}{2}}}d\zeta d\eta.
\eeq
In following subsections \ref{fh} and \ref{p-isv}, we use quaternionic analogue of Funk-Hecke formula to prove that
extremizer satisfying (\ref{sv-}) can only be constant (Theorem \ref{t-isv}). Concerning this involves most complicated computation, we leave it in independent subsections \ref{fh} and \ref{p-isv}.
Actually, we \textbf{claim} here: \\
\emph{
For any $4\le \lambda<Q=4n+6$ and any function $h$ that makes following integrals exist,
we have inverse inequality
\beq\label{isv}\iint_{S\times S}\frac{\overline{h(\zeta)}(\bar{\zeta}\cdot\eta+\bar{\eta}\cdot\zeta)h(\eta)}{|1-\zeta\cdot\bar{\eta}|^{\frac{\lambda}{2}}}d\zeta d\eta
\ge 2(p-1)\iint_{S\times S} \frac{\overline{h(\zeta)}h(\eta)}{|1-\zeta\cdot\bar{\eta}|^{\frac{\lambda}{2}}}d\zeta d\eta.\eeq
Moreover, ``=" holds if and only if
\[ h\in \left\{\begin{array}{ll}
        V_{0,0} & \lambda>4 \\
        V_{0,0}\bigoplus_{j\ge k \ge 2} V_{j,k} & \lambda=4.
      \end{array}\right.
\]
 }\\
If the \textbf{claim} is proved to be right, then for $4\le\lambda<Q$, above extremizer satisfying zero-center mass condition
can only be constant functions, i.e. $h\equiv const.$. Actually, for $\lambda=4$, from Theorem \ref{t-eig}, we know that
the eigenvalue $\lambda_{j,k}(|1-\zeta\cdot\bar{\eta}|^{-\frac{\lambda}{2}})$ vanishes for $k\ge 1$, so from Euler-Lagrange equation (\ref{EL}),
$h$ can only be constant. Another way to see this is considering the original functional of the inequality (\ref{shls})
as for nonnegative extremizer $h$ satisfying  zero center-mass condition, we have
\[\iint_{S\times S}\frac{h(\zeta)h(\eta)}{d_S(\zeta,\eta)^\lambda}d\zeta d\eta
=\iint_{S\times S}\frac{h_0(\zeta)h_0(\eta)}{d_S(\zeta,\eta)^\lambda}d\zeta d\eta,~\qquad \text{but}~\qquad\|h\|_p^2\ge\|h_0\|_p^2,\]
where $h_0=\frac{1}{|S|}\int_S h$ (mean of $h$ on $S$) is the projection of $h$ onto $V_{0,0}$.

\subsection{Sharp Constants and Extremizers}\label{sc}
Assume the \textbf{claim} holds, then we can give sharp constants and all the extremizers. \\
(1) we get the \emph{sharp constants} for (\ref{hls}) and (\ref{shls}).\\
Using the fact that $\int_S \frac{1}{|1-\eta_{n+1}|^{\frac{\lambda}{2}}}d\eta=\lambda_{0,0}(K_1^{\frac{\lambda}{4}})$
from (\ref{e-1}) in Theorem \ref{t-eig} and inserting extremizer $f=g=1$, we get the sharp constant
\begin{align*}
C'_\lambda=& 2^{\frac{\lambda}{2}}|S|^{1-\frac{2}{p}}\int_S \frac{1}{|1-\eta_{n+1}|^{\frac{\lambda}{2}}}d\eta\\
=&2^{\frac{\lambda}{2}}|S|^{1-\frac{2}{p}}
\frac{2\pi^{2n+2}\Gamma(2n-\frac{\lambda}{2}+3)}{\Gamma(2n-\frac{\lambda}{4}+2)\Gamma(2n-\frac{\lambda}{4}+3)},\\
C_\lambda=&2^{-\frac{(4n+3)\lambda}{Q}}C'_\lambda\\
=&2^{-\frac{2n\lambda}{Q}}|S|^{1-\frac{2}{p}}
\frac{2\pi^{2n+2}\Gamma(2n-\frac{\lambda}{2}+3)}{\Gamma(2n-\frac{\lambda}{4}+2)\Gamma(2n-\frac{\lambda}{4}+3)}.
\end{align*}
(2) we get all the \emph{extremizers} for  (\ref{hls}) and (\ref{shls}).\\
(I)\emph{extremizer for} (\ref{shls}).\\
On the one hand, for any extremizer $h$ of (\ref{shls}), $\tilde{h}$ is constant, which means $
h=|J_\gamma|^{\frac{1}{p}}$ is the only possible form of extremizer. To compute, first note
$|J_\gamma|=|J_{\ca{C}^{-1}}(A\zeta)||J_{\ca{C}}(A\gamma(\zeta))|, ~A\gamma(\zeta)=\ca{C}\circ \ca{S}_\delta \circ \ca{C}^{-1}(A\zeta)$, then
\begin{align*}
h ~\sim~ & \left(\frac{(1+|\delta q|^2)^2+|\delta^2 w|^2}{(1+|q|^2)^2+|w|^2}\right)^{-\frac{2Q-\lambda}{4}}\\
~=~ & \left|\frac{1+|\delta q|^2\pm\delta^2w}{1+|q|^2\pm w}\right|^{-\frac{2Q-\lambda}{2}}\\
~\sim~ & \left|1-\frac{\delta^2-1}{\delta^2+1}\frac{1-|q|^2\pm w}{1+|q|^2\mp w}\right|^{-\frac{2Q-\lambda}{2}} \\
~=~ & |1-\xi\cdot\bar{\zeta}|^{-\frac{2Q-\lambda}{2}},
\end{align*} with $(q,w)=\ca{C}^{-1}(A\zeta),~ \xi=A^{-1}(0,\ldots,\frac{\delta^2-1}{\delta^2+1})$, satisfying $|\xi|<1$.\\
On the other hand, given any $|\xi|<1,$ we can inverse the above process: there exists one $A\in Sp(n+1)$ s.t.
$A\xi=|\xi|(0,\ldots,0,1)$, so, through ``boundary" Cayley transform, we have
\begin{align*}
|J_{\ca{C}}|^{\frac{1}{p}}|1-\xi\cdot\overline{(A^{-1}\zeta)}|^{-\frac{2Q-\lambda}{2}}
~=~  & |J_{\ca{C}}|^{\frac{1}{p}}\left|1-|\xi|\frac{1-|q|^2\pm w}{1+|q|^2\mp w}\right|^{-\frac{2Q-\lambda}{2}}\\
~\sim~ & ((1+|q|^2)^2+|w|^2)^{-\frac{2Q-\lambda}{4}}\left(\frac{(1+|\delta q|^2)^2+|\delta^2 w|^2}{(1+|q|^2)^2+|w|^2}\right)^{-\frac{2Q-\lambda}{4}}\\
~=~ & ((1+|\delta q|^2)^2+|\delta^2 w|^2)^{-\frac{2Q-\lambda}{4}},
\end{align*}with $(q,w)= \ca{C}^{-1}\zeta,~\delta^2=\frac{1\pm|\xi|}{1\mp|\xi|}$, which tells that correspondence of function
$|1-\xi\cdot\overline{ (A^{-1}\zeta)}|^{-\frac{2Q-\lambda}{2}}$ through relation formula (\ref{corr}) is just constant multiple of $\delta-$dilation of particular extremizer $((1+|q|^2)^2+|w|^2)^{-\frac{2Q-\lambda}{4}}$ for (\ref{hls}), i.e. the function $|1-\xi\cdot\overline{ (A^{-1}\zeta)}|^{-\frac{2Q-\lambda}{2}}$ is an extremizer for (\ref{shls}). Because of the rotation-invariant of (\ref{shls}), any function $|1-\xi\cdot\bar{\zeta}|^{-\frac{2Q-\lambda}{2}}$, for any $|\xi|<1$ is an extremizer for (\ref{shls}).\\
Then, combine the two direction arguments, we claim: functions $\sim|1-\xi\cdot\bar{\zeta}|^{-\frac{2Q-\lambda}{2}}$, with $|\xi|<1$, are right all the extremizers for (\ref{shls}).\\
(II) \emph{extremizers for} (\ref{hls}).\\
For any $|\xi|<1$,
\begin{align*}
|J_{\ca{C}}|^{\frac{1}{p}}|1-\xi\cdot\bar{\zeta}|^{-\frac{2Q-\lambda}{2}}~\sim~ & ((1+|q|^2)^2+|w|^2)^{-\frac{2Q-\lambda}{4}}\left|1-\xi\cdot\overline{(\frac{2q}{1+|q|^2-w},\frac{1-|q|^2+w}{1+|q|^2-w})}\right|^{-\frac{2Q-\lambda}{2}}\\
~\sim~ & \left|1+|q|^2+w-(\xi_1,\xi_2,\ldots,\xi_n)\cdot\overline{(2q)}-\xi_{n+1}\cdot(1-|q|^2-w)\right|^{-\frac{2Q-\lambda}{2}}\\
~\sim~ & \left||q|^2+w-2q_0\cdot \bar{q}+r_0\right|^{-\frac{2Q-\lambda}{2}},
\end{align*} with $q_0=\frac{(\xi_1,\ldots,\xi_n)}{1+\xi_{n+1}}, r_0=\frac{1-\xi_{n+1}}{1+\xi_{n+1}}$, satisfying $\text{Re}r_0>|q_0|^2$. Then, because of the bijection between $\{\xi| |\xi|<1\}$ and $\{(q_0,p_0)| \text{Re}p_0>|q_0|^2\}$, we claim: all extremizers of (\ref{hls}) are explicitly given by
\beq\label{6} \sim \left||q|^2+w-2q_0\cdot \bar{q}+r_0\right|^{-\frac{2Q-\lambda}{2}},\eeq
 with $q_0\in\bb{H}^n, r_0\in\bb{H}$, satisfying $\text{Re}r_0>|q_0|^2$.
We remark that all the extremziers given by (\ref{6}) are just $c-$constant multiple, $\delta-$dilation and $u_0-$left translation
of the special extremizer $H(u)=\left((1+|q|^2)^2+|w|^2\right)^{-\frac{2Q-\lambda}{4}}$, which corresponds to constant function
extremizer for sphere inequality (\ref{shls}), with $\delta=(\text{Re}r_0-|q_0|^2)^{-\frac{1}{2}}, u_0=(q_0, -\text{Im}r_0)$.
Actually, for any $\delta>0, u_0=(q_0, w_0)\in G$, \begin{align*}
H(\delta( u_0^{-1}u))&=\left((1+|\delta(q-q_0)|^2)^2+|\delta^2(w-w_0-2\text{Im}q_0\cdot\bar{q})|^2\right)^{-\frac{2Q-\lambda}{4}}\\
&=|1+|\delta(q-q_0)|^2+\delta^2(w-w_0-2\text{Im}q_0\cdot\bar{q})|^{-\frac{2Q-\lambda}{2}}\\
&\sim||q|^2+w-2q_0\cdot\bar{q}+\frac{1}{\delta^2}+|q_0|^2-w_0|^{-\frac{2Q-\lambda}{2}}.
\end{align*} This gives the Remark (3) of Theorem \ref{t-q}.

The two main theorems are then proved after we check the \textbf{claim} in subsection \ref{ss-isv} which tells by
the following argument that extremizer that satifying zero center-mass condition $h$ can only be constant.
The \textbf{claim} is a corollary of Theorem \ref{t-isv} with $\alpha=\frac{\lambda}{4}$.

\subsection{Quaternionic Funk-Hecke Formula and Eigenvalues}\label{fh}
In this and the next subsection, we focus on proving the following bilinear estimate about inverse second-variation inequality.
\bthm\label{t-isv}\emph{[Bilinear Inequality]}\\
Let $1\le \alpha<\frac{Q}{4}$, then for any $f$ on $S$, we have
\beq\label{isv-}
\iint_{S\times S}\frac{\overline{f(\zeta)}(\bar{\zeta}\cdot\eta+\bar{\eta}\cdot\zeta)f(\eta)}{|1-\zeta\cdot\bar{\eta}|^{2\alpha}}d\zeta \eta\ge \frac{2\alpha}{\frac{Q}{2}-\alpha}\iint_{S\times S} \frac{\overline{f(\zeta)}f(\eta)}{|1-\zeta\cdot\bar{\eta}|^{2\alpha}}d\zeta d\eta,
\eeq and when $\alpha>1$, ``=" holds if and only if $f$ is constant function; When $\alpha=1$, "=" holds if and only if $f\in V_{0,0}\bigoplus_{j\ge k\ge 2}V_{j,k}$,
moreover, if $f$ is an extremizer for sharp HLS inequality \emph{(\ref{shls})}, $f$ can only be constant function.
Here, $\alpha\ge 1$ is sharp.\ethm
It suffices to prove the part before ``moreover" from the remark following the \textbf{claim} in subsection \ref{ss-isv}. First, in this subsection, we are going to derive an quaternionic analogue of classical real or complex Funk-Hecke theorem, which concerns the integral operators associated with kernel of the form $K(\zeta\cdot\bar{\eta})$. We then give eigenvalues of integral operators with two useful kernels. we left the prove of Theorem \ref{t-isv} to next subsection.

We have $O(4n+4)-$irreducible decomposition \[ L^2(S)=\bigoplus_{k\ge 0} \ca{H}_k,\]
with $\ca{H}_k$ the space of $k-$homogeneous harmonic polynomials in real variables,
and $U(2n+2)$-irreducible decomposition \[ L^2(S)=\bigoplus_{j,k\ge 0}\ca{H}_{j,k},\]
with $\ca{H}_{j,k}$ the space of $(j,k)-$bihomogeneous harmonic polynomials in complex variables and their complex conjugates.
For quaternionic case, we have $Sp(n+1)Sp(1)-$irreducible decomposition \beq\label{dec} L^2(S)=\bigoplus_{j\ge k\ge 0} V_{j,k},\eeq
where $V_{j,k}\subset \ca{H}_{j+k}$ is called ``$(j,k)-$bispherical harmonic space", which is generated from the action of $Sp(n+1)\times Sp(1)$ on zonal harmonic polynomial (Theorem 3.1 (4) in \cite{jw}, note that there is some print error in the formula, but one can easily correct from the proof)
\begin{align}\label{zon} Z_{j,k}(\zeta)~\sim~&\mbox{Proj}_{\ca{H}}\left(\sum_{l=k}^{[\frac{j+k}{2}]} (-1)^l C_{j+k-l}^l
 (2\mbox{Re} \zeta_{n+1})^{j+k-2l}|\zeta_{n+1}|^{2l}\right)\nonumber\\
~\sim~&\frac{k!(2n-1)!}{(2n+k-1)!}\frac{\sin(j-k+1)\phi}{\sin\phi}\cos^{j-k}\theta P_k^{(2n-1,j-k+1)}(\cos{2\theta}),
\end{align} with $|\zeta_{n+1}|=\cos\theta, \mbox{Re}\zeta_{n+1}=\cos\theta\cos\phi ~(\theta\in [0,\frac{\pi}{2}], \phi\in[0,\pi])$, where $\mbox{Proj}_{\ca{H}}$ is the projection to subspace of harmonic polynomials, $C_{j+k-l}^l$ is the combinatorial number and
$P_k^{(2n-1,j-k+1)}(z)$ is the Jacobi polynomial of order $k$ associated to weight $(1-z)^{2n-1}(1+z)^{j-k+1}$.
Here $``\sim"$ means equality modulo a constant multiple. Note that
$\sum_{l=0}^{[\frac{j+k}{2}]} (-1)^l C_{j+k-l}^l (2\mbox{Re} \zeta_{n+1})^{j+k-2l}|\zeta_{n+1}|^{2l} \in \ca{H}_{j+k}$.
See for example \cite{kostant1969existence,jw} for the classical spherical harmonic realization of spherical principle series.
\bthm\label{t-fh}\emph{[Quaternionic Funk-Hecke Formula]}\\
 Let $K$ be a function on unit ball of $\bb{H}$, s.t. the following integral exists, like $K\in L^1(B(0,1))$. Then integral operator with kernel of $K(\zeta\cdot\bar{\eta})$ is diagonal w.r.t decomposition \emph{(\ref{dec})}, and the eigenvalue on (j,k)-subspace $V_{j,k}$ is given by
\begin{align}
&~\lambda_{j,k}(K)\nonumber\\
=&~ \frac{2\pi^{2n}k!}{(j-k+1)(k+2n-1)!}\times \nonumber\\
&~ \int_0^{\frac{\pi}{2}} d\theta \sin^{4n-1}\theta\cos^{j-k+3}\theta P_k^{(2n-1,j-k+1)}(\cos{2\theta})\int_{\{u\in\bb{H}||u|=1\}\simeq \bb{S}^3} du K(\cos\theta u)\frac{\sin(j-k+1)\phi}{\sin\phi},
\end{align}
with $\re u=\cos\phi ~(\phi\in[0,\pi])$, $du$ is the standard Lebesgue measure on $\bb{S}^3$, i.e., $du=\sin^2\phi\sin\phi_2 d\phi d\phi_2 d\phi_1$,
in polar coordinates: $u=x_1+x_2i+(x_3+x_4i)j$, where
\begin{eqnarray*}
                                                                      x_4 &=& \sin\phi\sin\phi_2\sin\phi_1 \\
                                                                      x_3 &=& \sin\phi\sin\phi_2\cos\phi_1 \\
                                                                      x_2 &=& \sin\phi\cos\phi_2 \\
                                                                      x_1 &=& \cos\phi \qquad (\phi,\phi_2\in[0,\pi],\phi_3\in[0,2\pi]).
                                                                    \end{eqnarray*}

\ethm
\bpf
From Schur's lemma and the irreducibility of $(j,k)$-subspace $V_{j,k}$, we see the integral operator associated to
$K(\zeta\cdot\bar{\eta})$ is diagonal with eigenvalues denoted by $\lambda_{j,k}$. Now, we compute the eigenvalues.
Assume $\{Y_{j,k}^\mu\}_{1\le\mu\le m_{j,k}}$ is a normalized orthogonal basis of $V_{j,k}$, then in abuse of notation
the reproducing kernel of projection operator onto $V_{j,k}$ is given by
\[Z_{j,k}(\zeta,\eta)=Z_{j,k}(\zeta\cdot\bar{\eta})=\sum_{\mu=1}^{m_{j,k}}Y_{j,k}^\mu(\zeta) \overline{Y_{j,k}^\mu(\eta)}.\]
Then we have
\[\int_S K(\zeta\cdot\bar{\eta})Z_{j,k}(\eta\cdot\bar{\zeta})d\eta=\lambda_{j,k}Z_{j,k}(1),\] which implies
\begin{align}\label{eig}
\lambda_{j,k} &= Z^{-1}_{j,k}(1)\int_S K(\zeta\cdot\bar{\eta})Z_{j,k}(\eta\cdot\bar{\zeta})d\eta\nonumber\\
&=Z^{-1}_{j,k}(1)\int_S K(\overline{\eta_{n+1}})Z_{j,k}(\eta_{n+1})d\eta.
\end{align}
In polar coordinates for $\eta=(\eta_1,\eta_2,\ldots,\eta_{n+1})$,
\begin{eqnarray*}
  \eta_1 &=& u_1\sin\theta_n\sin\theta_{n-1}\ldots\sin\theta_1 \\
  \eta_2 &=& u_2\sin\theta_n\sin\theta_{n-1}\ldots\cos\theta_1\\
  \eta_i &=& \ldots \\
  \eta_n &=& u_n\sin\theta_n\cos\theta_{n-1} \\
  \eta_{n+1} &=& u_{n+1}\cos\theta_n,\\
    u_i^4 &=& \sin\phi_i^3\sin\phi_i^2\sin\phi_i^1 \\
    u_i^3 &=& \sin\phi_i^3\sin\phi_i^2\cos\phi_i^1 \\
    u_i^2 &=& \sin\phi_i^3\cos\phi_i^2 \\
    u_i^1 &=& \cos\phi_i^3
\end{eqnarray*} with $\theta_i\in [0,\frac{\pi}{2}], \phi_i^3,\phi_i^2\in[0,\pi],\phi_i^1\in[0,2\pi], ~\forall 1\le i\le n+1$,
 we have the invariance measure (11.7.3 (2) in \cite{vk})
\beq\label{mea} d\eta= 
\prod_{i=1}^n (\sin^{4i-1}\theta_i\cos^3\theta_id\theta_i)\prod_{j=1}^{n+1}du_j,
\eeq with $du_j=\sin^2\phi_{j}^3\sin\phi_{j}^2d\phi_{j}^3d\phi_{j}^3d\phi_{j}^2d\phi_{j}^1$.
Putting the formula (\ref{zon}) for zonal harmonics and invariance measure (\ref{mea}) on $S$ into (\ref{eig}), and from $|\bb{S}^{4n-1}|=\frac{2\pi^{2n}}{(2n)!}$, $P_k^{(2n-1,j-k+1)}(1)=\frac{(k+2n-1)!}{k!(2n-1)!}$ (22.2.1 in \cite{abramowitz2012handbook}), we get
\begin{align*}
&~\lambda_{j,k} \\
=&~ |\bb{S}^{4n-1}|((j-k+1)P_k^{(2n-1,j-k+1)}(1))^{-1}\times\\
&~\int_0^{\frac{\pi}{2}} d\theta_n \sin^{4n-1}\theta_n\cos^{j-k+3}\theta_n P_k^{(2n-1,j-k+1)}(\cos{2\theta_n})\int_{\bb{S}^3} du_{n+1} K(\cos\theta_n u_{n+1})\frac{\sin(j-k+1)\phi_{n+1}^3}{\sin\phi_{n+1}^3}\\
=&~ \frac{2\pi^{2n}k!}{(j-k+1)(k+2n-1)!}\times \\
&~ \int_0^{\frac{\pi}{2}} d\theta_n \sin^{4n-1}\theta_n\cos^{j-k+3}\theta_n P_k^{(2n-1,j-k+1)}(\cos{2\theta_n})\int_{\bb{S}^3} du_{n+1} K(\cos\theta_n u_{n+1})\frac{\sin(j-k+1)\phi_{n+1}^3}{\sin\phi_{n+1}^3}
\end{align*}
\epf

In order to prove Theorem \ref{t-isv}, it suffices to compute eigenvalue of functions of two forms
$K^\alpha_1(q)=|1-q|^{-2\alpha},K^\alpha_2(q)=|q|^2|1-q|^{-2\alpha}$, noting that
$\zeta\cdot\bar{\eta}+\eta\cdot\bar{\zeta}=2\re\zeta\cdot\bar{\eta}=1+|\zeta\cdot\bar{\eta}|^2-|1-\zeta\cdot\bar{\eta}|^2$.
\bthm\emph{[Eigenvalues]}\label{t-eig}\\
Given $-1<\alpha<\frac{Q}{4}$. Denote \emph{(15.1.1 and 15.1.20 in \cite{abramowitz2012handbook})}
\begin{align}\label{A}
A(a,b,c) &= \sum_{\mu\ge 0}\frac{\Gamma(\mu+a)\Gamma(\mu+b)}{\mu!\Gamma(\mu+c)}=\frac{\Gamma(a)\Gamma(b)\Gamma(c-a-b)}{\Gamma(c-a)\Gamma(c-b)}, ~c>a+b; \nonumber\\
(a,b,c) &= (j+\alpha,k+\alpha-1,j+k+\frac{Q}{2}-1).
\end{align} Then\\
\emph{(1)} The eigenvalues of integral operators associated to kernel $K^\alpha_1(q)=|1-q|^{-2\alpha}$ are given by
\begin{align}\label{e-1}
\lambda_{j,k}(K^\alpha_1)&=~\frac{2\pi^{2n+2}}{\Gamma(\alpha)\Gamma(\alpha-1)} A(a,b,c)\nonumber\\
&=~2\pi^{2n+2}\frac{\Gamma(\frac{Q}{2}-2\alpha)}{\Gamma(\alpha)\Gamma(\alpha-1)}
\frac{\Gamma(j+\alpha)}{\Gamma(j+\frac{Q}{2}-\alpha)}
\frac{\Gamma(k+\alpha-1)}{\Gamma(k+\frac{Q}{2}-\alpha-1)}.
\end{align}
\emph{(2)} The eigenvalues of integral operators associated to kernel $K^\alpha_2(q)=|q|^2|1-q|^{-2\alpha}$ are given by
\beq\lambda_{j,k}(K^\alpha_2)=C_{j,k}^{\alpha}\lambda_{j,k}(K^\alpha_1),\eeq with
\begin{align} &C_{j,k}^{\alpha}\nonumber\\
=&1-(\alpha-2)(c-a-b)\left(\frac{1}{(a-1)(c-a)}+\frac{1}{(b-1)(c-b)}-(\alpha-2)\frac{c-a-b+1}{(a-1)(b-1)(c-a)(c-b)}\right)\nonumber\\
=&1-(\alpha-2)(\frac{Q}{2}-2\alpha)\times \nonumber\\
&\frac{-(j+\alpha)^2-(k+\alpha-1)^2+(j+k+\frac{Q}{2})(j+2\alpha+k-1)-2(j+k+\frac{Q}{2}-1)-(\alpha-2)(\frac{Q}{2}-2\alpha+1)}
{(j+\alpha-1)(k+\alpha-2)(k+\frac{Q}{2}-\alpha-1)(j+\frac{Q}{2}-\alpha)}.
\end{align}  In the singular point $\alpha=0,1,2$, the above formula can be viewed as limit, fixing $j,k$.
\ethm
\bpf
(1)Putting $K=K^\alpha_1$ into Theorem \ref{fh}, we get
\begin{align}\label{1}
\lambda_{j,k}(K^\alpha_1)=&~\frac{4\pi^{2n+1}k!}{(j-k+1)(k+2n-1)!}\int_0^{\frac{\pi}{2}} d\theta \sin^{4n-1}\theta\cos^{j-k+3}\theta P_k^{(2n-1,j-k+1)}(\cos{2\theta})\times\nonumber\\
&~\int_0^\pi d\phi (1+\cos^2\theta-2\cos\phi\cos\theta)^{-\alpha}(\cos(j-k)\phi-\cos(j-k+2)\phi).
\end{align}
Using Gegenbauer polynomials, the following fact exists ((5.11) in \cite{fl}, \cite{abramowitz2012handbook})
\begin{align}\label{int}&~\int_0^\pi d\phi (1+\cos^2\theta-2\cos\phi\cos\theta)^{-\alpha}\cos(j-k)\phi\nonumber\\
=&~\frac{\pi}{\Gamma^2(\alpha)}\sum_{\mu\ge 0}\cos^{|j-k|+2\mu}\theta\frac{\Gamma(\mu+\alpha)\Gamma(\mu+|j-k|+\alpha)}{\mu!(\mu+|j-k|)!},\end{align} which, from (\ref{1}) and $j\ge k\ge 0$, gives
\begin{align}\label{lam}
\lambda_{j,k}(K^\alpha_1)=&~\frac{4\pi^{2n+2}k!}{(j-k+1)(k+2n-1)!\Gamma^2(\alpha)}\times \nonumber\\
&~\bigg(\sum_{\mu\ge 0}\frac{\Gamma(\mu+\alpha)\Gamma(\mu+|j-k|+\alpha)}{\mu!(\mu+|j-k|)!}\int_0^{\frac{\pi}{2}} d\theta \sin^{4n-1}\theta\cos^{2(j-k)+3+2\mu}\theta P_k^{(2n-1,j-k+1)}(\cos{2\theta})- \nonumber\\
&~\sum_{\mu\ge 0}\frac{\Gamma(\mu+\alpha)\Gamma(\mu+|j-k|+2+\alpha)}{\mu!(\mu+|j-k|+2)!}\int_0^{\frac{\pi}{2}} d\theta \sin^{4n-1}\theta\cos^{2(j-k)+5+2\mu}\theta P_k^{(2n-1,j-k+1)}(\cos{2\theta})\bigg) \nonumber\\
\triangleq&~\frac{4\pi^{2n+2}k!}{(j-k+1)(k+2n-1)!\Gamma^2(\alpha)}(I_1-I_2).
\end{align}
Using the following Rodrigues' formula ((22.11.1) in \cite{abramowitz2012handbook}),
\[P_k^{(2n-1,j-k+1)}(t)=\frac{(-1)^k}{2^kk!}(1-t)^{-(2n-1)}(1+t)^{-(j-k+1)}\frac{d^k}{dt^k}\{(1-t)^{2n+k-1}(1+t)^{j+1}\},\]
changing variable $\cos2\theta=t$ and integrating by part, we get
\begin{align}\label{lam-}
&~\int_0^{\frac{\pi}{2}} d\theta \sin^{4n-1}\theta\cos^{2(j-k)+3+2\mu}\theta P_k^{(2n-1,j-k+1)}(\cos{2\theta})\nonumber\\
=&~\frac{(-1)^k2^{-(\mu+j+2n+2)}}{k!}\int_{-1}^1 dt (1+t)^\mu\frac{d^k}{dt^k}\{(1-t)^{2n+k-1}(1+t)^{j+1}\}\nonumber\\
=&~\chi_{\mu\ge k}\frac{2^{-(\mu+j+2n+2)}\mu!}{k!(\mu-k)!}\int_{-1}^1 dt (1+t)^{\mu+j-k+1}(1-t)^{2n+k-1}\nonumber\\
=&~\chi_{\mu\ge k}2^{-1}\frac{\mu!B(\mu+j-k+2,2n+k)}{k!(\mu-k)!}\nonumber\\
=&~\chi_{\mu\ge k}2^{-1}\frac{\mu!\Gamma(\mu+j-k+2)\Gamma(2n+k)}{k!(\mu-k)!\Gamma(\mu+j+2n+2)}.
\end{align}
Inserting (\ref{lam-}) into (\ref{lam}), we get
\begin{align}\label{I}
I_1=&~2^{-1}\sum_{\mu\ge k}
\frac{\Gamma(\mu+\alpha)\Gamma(\mu+|j-k|+\alpha)}{\mu!(\mu+|j-k|)!}
\frac{\mu!\Gamma(\mu+j-k+2)\Gamma(2n+k)}{k!(\mu-k)!\Gamma(\mu+j+2n+2)}\nonumber\\
I_2=&~2^{-1}\sum_{\mu+1\ge max\{k,1\}}
\frac{\Gamma(\mu+\alpha)\Gamma(\mu+|j-k|+2+\alpha)}{\mu!(\mu+|j-k|+2)!}
\frac{(\mu+1)!\Gamma(\mu+1+j-k+2)\Gamma(2n+k)}{k!(\mu+1-k)!\Gamma(\mu+1+j+2n+2)}.
\end{align}Now, we get the differnce of the two terms.\\
For $k\ge 1$,
\begin{align*}
I_1-I_2=&~2^{-1}\sum_{\mu\ge k}
\frac{\mu!\Gamma(\mu+j-k+2)\Gamma(2n+k)}{k!(\mu-k)!\Gamma(\mu+j+2n+2)}\times\nonumber\\
&~\left(\frac{\Gamma(\mu+\alpha)\Gamma(\mu+|j-k|+\alpha)}{\mu!(\mu+|j-k|)!}-
\frac{\Gamma(\mu-1+\alpha)\Gamma(\mu+|j-k|+1+\alpha)}{(\mu-1)!(\mu+|j-k|+1)!}\right)\\
=&~2^{-1}\sum_{\mu\ge k}
\frac{\mu!\Gamma(\mu+j-k+2)\Gamma(2n+k)}{k!(\mu-k)!\Gamma(\mu+j+2n+2)}\frac{(\alpha-1)(j-k+1)\Gamma(\mu+\alpha-1)\Gamma(\mu+j-k+\alpha)}{\mu!(\mu+j-k+1)!}\\
=&~\frac{2^{-1}(\alpha-1)(j-k+1)\Gamma(2n+k)}{k!}\sum_{\mu\ge k}\frac{\Gamma(\mu+\alpha-1)\Gamma(\mu+j-k+\alpha)}{(\mu-k)!\Gamma(\mu+j+2n+2)}\nonumber\\
=&~\frac{2^{-1}(\alpha-1)(j-k+1)\Gamma(2n+k)}{k!}\sum_{\mu\ge 0}\frac{\Gamma(\mu+k+\alpha-1)\Gamma(\mu+j+\alpha)}{\mu!\Gamma(\mu+j+k+2n+2)};
\end{align*}
For $k=0$,
\begin{align*}
I_1-I_2=&~2^{-1}(\alpha-1)(j+1)\Gamma(2n)\sum_{\mu\ge 1}\frac{\Gamma(\mu+\alpha-1)\Gamma(\mu+j+\alpha)}{\mu!\Gamma(\mu+j+2n+2)}
+2^{-1}\Gamma(2n)\frac{\Gamma(\alpha)\Gamma(j+\alpha)(j+1)}{\Gamma(j+2n+2)}\\
=&~2^{-1}(\alpha-1)(j+1)\Gamma(2n)\sum_{\mu\ge 0}\frac{\Gamma(\mu+\alpha-1)\Gamma(\mu+j+\alpha)}{\mu!\Gamma(\mu+j+2n+2)}
\end{align*}
So, together with (\ref{A})and (\ref{lam}), we get\\
\begin{align*}
\lambda_{j,k}(K_1^\alpha)=&~\frac{2\pi^{2n+2}}{\Gamma(\alpha)\Gamma(\alpha-1)} A(a,b,c)\\
=&~2\pi^{2n+2}\frac{\Gamma(\frac{Q}{2}-2\alpha)}{\Gamma(\alpha)\Gamma(\alpha-1)}
\frac{\Gamma(j+\alpha)}{\Gamma(j+\frac{Q}{2}-\alpha)}
\frac{\Gamma(k+\alpha-1)}{\Gamma(k+\frac{Q}{2}-\alpha-1)}.
\end{align*}
The fist part of theorem is proved.\\
(2) Putting $K=K^\alpha_2$ into Theorem \ref{t-fh}, we get
\begin{align}
\lambda_{j,k}(K^\alpha_2)=&~\frac{4\pi^{2n+1}k!}{(j-k+1)(k+2n-1)!}\int_0^{\frac{\pi}{2}} d\theta \sin^{4n-1}\theta\cos^{j-k+5}\theta P_k^{(2n-1,j-k+1)}(\cos(2\theta))\times\nonumber\\
&~\int_0^\pi d\phi (1+\cos^2\theta-2\cos\phi\cos\theta)^{-\alpha}(\cos(j-k)\phi-\cos(j-k+2)\phi).
\end{align} Compare with (\ref{1}), $\mu$ can be substituted by $\mu+1$, so repeat the same computation in (1), we obtain analogue of (\ref{lam}), with
\begin{align}\label{I-}
I_1=&~2^{-1}\sum_{\mu+1\ge max\{k,1\}}
\frac{\Gamma(\mu+\alpha)\Gamma(\mu+|j-k|+\alpha)}{\mu!\Gamma(\mu+|j-k|)}
\frac{(\mu+1)!\Gamma(\mu+1+j-k+2)\Gamma(2n+k)}{k!(\mu+1-k)!\Gamma(\mu+1+j+2n+2)}\nonumber\\
I_2=&~2^{-1}\sum_{\mu+2\ge max\{k,2\}}
\frac{\Gamma(\mu+\alpha)\Gamma(\mu+|j-k|+2+\alpha)}{\mu!\Gamma(\mu+|j-k|+2)}
\frac{(\mu+2)!\Gamma(\mu+2+j-k+2)\Gamma(2n+k)}{k!(\mu+2-k)!\Gamma(\mu+2+j+2n+2)}.
\end{align}
For $k\ge 2$,
\begin{align}
\lambda_{j,k}(K^\alpha_2)=&~\frac{2\pi^{2n+2}}{\Gamma(\alpha)\Gamma(\alpha-1)}\sum_{\mu\ge 0}\frac{(\mu+k)(\mu+j+1)\Gamma(\mu+k+\alpha-2)\Gamma(\mu+j+\alpha-1)}{\mu!\Gamma(\mu+j+k+2n+2)}\nonumber\\
=&~\frac{2\pi^{2n+2}}{\Gamma(\alpha)\Gamma(\alpha-1)}\sum_{\mu\ge 0}\frac{\Gamma(\mu+k+\alpha-1)\Gamma(\mu+j+\alpha)}{\mu!\Gamma(\mu+j+k+2n+2)}\frac{(\mu+k)(\mu+j+1)}{(\mu+k+\alpha-2)(\mu+j+\alpha-1)}.
\end{align}
Noting that
\begin{align}
&~\frac{(\mu+k)(\mu+j+1)}{(\mu+k+\alpha-2)(\mu+j+\alpha-1)}=~\frac{(\mu+a-1-(\alpha-2))(\mu+b-1-(\alpha-2))}{(\mu+a-1)(\mu+b-1)}\nonumber\\
=&~1-(\alpha-2)\left(\frac{1}{\mu+a-1}+\frac{1}{\mu+b-1}-(\alpha-2)\frac{1}{(\mu+a-1)(\mu+b-1)}\right),
\end{align}
so, from (\ref{A}), we have
\beq\lambda_{j,k}(K^\alpha_2)=C_{j,k}^{\alpha}\lambda_{j,k}(K^\alpha_1),\eeq with
\beq C_{j,k}^{\alpha}=1-(\alpha-2)(c-a-b)\left(\frac{1}{(a-1)(c-a)}+\frac{1}{(b-1)(c-b)}-(\alpha-2)\frac{c-a-b+1}{(a-1)(b-1)(c-a)(c-b)}\right).
\eeq
For $k<2(k=0,1)$, we only prove for $k=1$ as they are similar, using the idea as in (1) that the special degenerate case can be formally integrated into the general case. We denote $(I_1-I_2)_i$ the obvious term for $K^\alpha_i$, see (\ref{I}) and (\ref{I-}), then
\begin{align*}
(I_1-I_2)_2=&~(I_1-I_2)_1C_{j,k}^{\alpha}+(\alpha-1)(j+1)(\alpha-2)\times\\
&~\left(\frac{\Gamma(\alpha-1)\Gamma(j+\alpha)}{\Gamma(j+2n+3)}
+\frac{\Gamma(\alpha)\Gamma(j+\alpha-1)}{\Gamma(j+2n+3)}
-(\alpha-2)\frac{\Gamma(\alpha-1)\Gamma(j+\alpha-1)}{\Gamma(j+2n+3)}\right)-\\
&~\frac{\Gamma(j+2)}{\Gamma(j+2n+3)}\left(
\frac{(\alpha-1)j\Gamma(\alpha)\Gamma(j+\alpha)}{(j+1)!}-
\frac{\Gamma(\alpha)\Gamma(j+\alpha-1)}{(j-1)!}\right)\\
=&~C_{j,k}^{\alpha}(I_1-I_2)_1.
\end{align*}
The second part of theorem is proved.
\epf

\subsection{Proof of Bilinear Inequality}\label{p-isv}
We now proceed the \emph{proof of Theorem} \ref{t-isv}:\\
\bpf
From Theorem \ref{t-eig}, we have
\[\lambda_{j,k}(K^{\alpha-1}_1)=\frac{2\pi^{2n+2}}{\Gamma(\alpha-1)\Gamma(\alpha-2)} A(a-1,b-1,c)
=\lambda_{j,k}(K^{\alpha}_1)(\alpha-1)(\alpha-2)\frac{(c-a-b)(c-a-b+1)}{(a-1)(b-1)(c-a)(c-b)}.\]
So, from $\zeta\cdot\bar{\eta}+\eta\cdot\bar{\zeta}=1+|\zeta\cdot\bar{\eta}|^2-|1-\zeta\cdot\bar{\eta}|^2$,
to prove the theorem, it suffices to check:

When $\alpha>1$,
\[2-(\alpha-2)(c-a-b)\left(\frac{1}{(a-1)(c-a)}+\frac{1}{(b-1)(c-b)}+\frac{c-a-b+1}{(a-1)(b-1)(c-a)(c-b)}\right)\ge \frac{2\alpha}{\frac{Q}{2}-\alpha},\] which is
\beq\label{2} (\alpha-2)\left(\frac{1}{(a-1)(c-a)}+\frac{1}{(b-1)(c-b)}+\frac{c-a-b+1}{(a-1)(b-1)(c-a)(c-b)}\right)\le \frac{2}{\frac{Q}{2}-\alpha},\eeq as $c-a-b=\frac{Q}{2}-2\alpha>0$.\\
For $\alpha-2> 0$, then it suffices to check
\[\frac{2}{(\alpha-2)(\frac{Q}{2}-\alpha)}-
\left(\frac{1}{(a-1)(c-a)}+\frac{1}{(b-1)(c-b)}\right)\ge \frac{c-a-b+1}{(a-1)(b-1)(c-a)(c-b)}\]
Substitute $(a,b,c)$ (see (\ref{A})), the inequality becomes
\begin{align}\label{3}
&\frac{j(\alpha-2)+k(\frac{Q}{2}-\alpha)+kj}{(\alpha-2)(\frac{Q}{2}-\alpha)(\alpha-2+k)(\frac{Q}{2}-\alpha+j)}
+\frac{(k-1)(\alpha-2)+(j+1)(\frac{Q}{2}-\alpha)+(j+1)(k-1)}{(\alpha-2)(\frac{Q}{2}-\alpha)(\alpha-2+j+1)(\frac{Q}{2}-\alpha+k-1)}\nonumber\\
&\ge \frac{\frac{Q}{2}-2\alpha+1}{(\alpha-2+k)(\frac{Q}{2}-\alpha+j)(\alpha-2+j+1)(\frac{Q}{2}-\alpha+k-1)}
\end{align}
First note that, in (\ref{3}), the denominator of left side is smaller than that in the right side,  precisely i.e.
\[(\alpha-2)(\frac{Q}{2}-\alpha) \le min\{ (\alpha-2+k)(\frac{Q}{2}-\alpha+j) ,(\alpha-2+j+1)(\frac{Q}{2}-\alpha+k-1) \},\] which is equivalent to \[(k-1)(\alpha-2)+(j+1)(\frac{Q}{2}-\alpha)+(j+1)(k-1)\ge 0,\] and we easily see that the left side $\ge \frac{Q}{2}-\alpha-1-(\alpha-2)=\frac{Q}{2}-2\alpha+1 \ge 1$.\\
Then compare the numerators in (\ref{3}), we have that the sum of numerators in the left side is larger than that in the right side, i.e.
\beq\label{4} (j+k-1)(\alpha-2)+(k+j+1)(\frac{Q}{2}-\alpha)+kj+(j+1)(k-1)\ge \frac{Q}{2}-2\alpha+1,\eeq
which is just \[(\frac{Q}{2}-3)j+(\frac{Q}{2}-1)k+2jk\ge 0.\] This inequality is obviously right and become equality if and only if $j=k=0$, and in this special case, inequalities above all and (\ref{3}) become equality. So, for $\alpha>2$, we have proved that the equality for (\ref{isv-}) holds if and only if $f$ is a constant.\\
For $\alpha-2<0$, from (\ref{2}), it suffices to prove the inverse of (\ref{3}),
\begin{align}
&\frac{j(\alpha-2)+k(\frac{Q}{2}-\alpha)+kj}{(\alpha-2)(\frac{Q}{2}-\alpha)(\alpha-2+k)(\frac{Q}{2}-\alpha+j)}
+\frac{(k-1)(\alpha-2)+(j+1)(\frac{Q}{2}-\alpha)+(j+1)(k-1)}{(\alpha-2)(\frac{Q}{2}-\alpha)(\alpha-2+j+1)(\frac{Q}{2}-\alpha+k-1)}\nonumber\\
&\le \frac{\frac{Q}{2}-2\alpha+1}{(\alpha-2+k)(\frac{Q}{2}-\alpha+j)(\alpha-2+j+1)(\frac{Q}{2}-\alpha+k-1)}.
\end{align}
For $k\ge 1$, it's obvious by checking sign that strict inequality holds.
For $k=0$, it's also obvious and ``=" holds if and only if $j=0$, checking sign of fist term and using (\ref{4}) for $j\ge 1$.
$\alpha=2$ case is obvious by checking similarly or limit argument.

When $\alpha\le 1$, for Theorem \ref{t-isv}, we need to check
\begin{align}\label{<1}
&\Bigg(\frac{j(\alpha-2)+k(\frac{Q}{2}-\alpha)+kj}{(\alpha-2)(\frac{Q}{2}-\alpha)(\alpha-2+k)(\frac{Q}{2}-\alpha+j)}
+\frac{(k-1)(\alpha-2)+(j+1)(\frac{Q}{2}-\alpha)+(j+1)(k-1)}{(\alpha-2)(\frac{Q}{2}-\alpha)(\alpha-2+j+1)(\frac{Q}{2}-\alpha+k-1)}-\nonumber\\
&\frac{\frac{Q}{2}-2\alpha+1}{(\alpha-2+k)(\frac{Q}{2}-\alpha+j)(\alpha-2+j+1)(\frac{Q}{2}-\alpha+k-1)}\Bigg)
\frac{\Gamma(k+\alpha-1)}{\Gamma(\alpha-1)}\le 0.
\end{align}
When $\alpha=1$, (\ref{<1}) holds and equality holds iff $k=j=0$ or $k\ge 2$.
When $\alpha<1$, we find it still holds for $k=0$ and equality holds iff $j=0$. For $k=1$, the opposite inequality holds strictly for
$j\ge\frac{\frac{Q}{2}-\alpha}{1-\alpha}$,
but uncertain for $j<\frac{\frac{Q}{2}-\alpha}{1-\alpha}$. Anyway, there exists a $j_\alpha$, s.t. (\ref{<1}) holds for $j\le j_\alpha$, especially holds strictly when $j=1$, but for  $j> j_\alpha$, the opposite inequality holds. For $k\ge 2$, the opposite inequality holds strictly.\\
Then Theorem \ref{t-isv} is proved for $\alpha\ge 1$, which is sharp range.
\epf
\textbf{Remark:} \emph{local extremizer for small $\lambda$.}
To figure out the extremizer for $\lambda<4$,
we want to find counterexample or add more condition to restrict the extremizer to be constant function. But what the extremizer could be is still unknown by us so far. However, we remark substitutively that constant function is indeed ``local" extremizer. More precisely, we proved that (\ref{sv}) holds strictly for $h=1$ and all $\varphi$ satisfying $\int_S \varphi=0, \int_S \varphi(\zeta)\zeta d\zeta=0$. Note the second assumption of $\varphi$ comes from the conformal symmetry group. If we decompose it w.r.t (\ref{dec}), $\varphi=\sum_{j,k} Y_{j,k}$ with $Y_{j,k} \in V_{j,k}, Y_{0,0}=Y_{1,0}=0$ ($V_{0,0}$ is constant function space, $V_{1,0}$ is spanned by $\re \zeta_j, 1\le j\le n+1$), it suffices to prove
\[\sum_{j,k}(\lambda_{j,k}^\alpha-(p-1)\lambda_{0,0}^\alpha)|S|\norm{Y_{j,k}}_2^2<0,\] where $\lambda_{j,k}^\alpha$ is the eigenvalue of operator $|1-\zeta\cdot\overline{\eta}|^{-2\alpha}$ and $\alpha=\frac{\lambda}{4}$.
From (\ref{e-1}) in Theorem \ref{t-eig} (1), we can check that
\[\lambda_{j,k}^\alpha-(p-1)\lambda_{0,0}^\alpha\] is positive when $j=k=0$, zero when $j=1,k=0$, and negative otherwise. So, we have proved the negative second variation and therefore ``local" maximum at constant function 1.

\section*{Acknowledgement}
The work was partially done during An Zhang's visit as a joint Ph.D. student in the department of mathematics, UC, Berkeley. He would like to thank especially the department for hospitality.

\bibliographystyle{amsalpha}
\bibliography{hls}

\vskip 3\baselineskip
\flushleft

Michael Christ \\
Department of Mathematics,
University of California, Berkeley;
Berkeley, CA, USA.\\
\emph{email:} \texttt{mchrist@math.berkeley.edu}

\vskip 2\baselineskip

Heping Liu\\
School of Mathematical Science,
Peking University;
Beijing, China.\\
\emph{email:} \texttt{hpliu@math.pku.edu.cn}

\vskip 2\baselineskip

An Zhang\\
School of Mathematical Science,
Peking University;
Beijing, China.\\
\emph{email:} \texttt{anzhang@pku.edu.cn}\\
Current Address:\\
Department of Mathematics,
University of California, Berkeley;
Berkeley, CA, USA. \\
\emph{email:} \texttt{anzhang@math.berkeley.edu}

\end{document}